\documentclass[a4paper,fleqn,leqno]{article}

\usepackage{latexsym,amssymb,amsthm,amsmath,amsfonts,mathtools,epsfig,color,epsf,graphicx}
\usepackage[english]{babel}

\usepackage[utf8]{inputenc} 
\usepackage[T1]{fontenc}  
\usepackage{lmodern}  
\usepackage{tikz}
\usetikzlibrary{trees}
\usetikzlibrary{decorations.pathreplacing}
\usepackage{tabu}
\usepackage{array,multirow,multicol}
\usepackage{hyperref}   
\usepackage{enumerate}    
\usepackage{url}            
\usepackage{booktabs}       
 \usepackage{nicefrac}       
\usepackage{microtype}      
\usepackage{lipsum}

\newcommand{\N}{\mathbb{N}}
\newcommand{\R}{\mathbb{R}}
\renewcommand{\P}{\mathbb{P}}
\newcommand{\Z}{\mathbb{Z}}
\newcommand{\E}{\mathbb{E}}

\newcommand{\X}{\mathcal{X}}
\renewcommand{\S}{\mathcal{S}}

\renewcommand{\d}{{\rm d}}
\newcommand{\e}{{\rm e}}
\newcommand{\ii}{{\rm i}}

\makeatletter
 \@addtoreset{equation}{section}
 \makeatother

 \newcounter{extralabel}[section]

 \newtheorem{ittheorem}{Theorem}
 \newtheorem{itlemma}{Lemma}
 \newtheorem{itproposition}{Proposition}
 \newtheorem{itdefinition}{Definition}
 \newtheorem{itcorollary}{Corollary}
 \newtheorem{itconjecture}{Conjecture}
 \newtheorem{itremark}{Remark}
 \theoremstyle{definition}
 \newtheorem{itassumption}{Assumption}

 \newenvironment{theorem}{\addtocounter{extralabel}{1}
 \begin{ittheorem}}{\end{ittheorem}}

 \newenvironment{lemma}{\addtocounter{extralabel}{1}
 \begin{itlemma}}{\end{itlemma}}

 \newenvironment{proposition}{\addtocounter{extralabel}{1}
 \begin{itproposition}}{\end{itproposition}}

 \newenvironment{definition}{\addtocounter{extralabel}{1}
 \begin{itdefinition}}{\end{itdefinition}}

 \newenvironment{corollary}{\addtocounter{extralabel}{1}
 \begin{itcorollary}}{\end{itcorollary}}

 \newenvironment{remark}{\addtocounter{extralabel}{1}
 \begin{itremark}}{\end{itremark}}

\newenvironment{assumption}{\addtocounter{extralabel}{1}
	\begin{itassumption}}{\end{itassumption}}

\tikzstyle{level 1}=[level distance=6em, sibling distance=10em]
\tikzstyle{level 2}=[level distance=3.5cm, sibling distance=2cm]
\tikzstyle{bag} = [text width=4em, text centered]
\tikzstyle{end} = [circle, minimum width=3pt,fill, inner sep=0pt]
\tikzset{dummy/.style = {circle,draw}}

\title{Spatially inhomogeneous populations with seed-banks:\\ 
II. Clustering regime}

\begin{document}
\author{Frank den Hollander$^1$, Shubhamoy Nandan$^2$}
\date{\today}

\maketitle

\begin{abstract}
We consider a spatial version of the classical Moran model with seed-banks where the constituent populations have finite sizes. Individuals live in colonies labelled by $\Z^d$, $d \geq 1$, playing the role of a geographic space, {carry one of two \emph{types}, $\heartsuit$ or $\spadesuit$,}  and change type via \emph{resampling} as long as they are \emph{active}. Each colony contains a seed-bank into which individuals can enter to become \emph{dormant}, suspending their resampling until they exit the seed-bank and become active again. Individuals resample not only from their own colony, but also from other colonies according to a symmetric random walk transition kernel. The latter is referred to as \emph{migration}. The sizes of the active and the dormant populations depend on the colony and remain constant throughout the evolution.

It was shown in \cite{HN01} that the spatial system is well-defined, {admits a family of equilibria parametrized by the initial density of type $\heartsuit$}, and exhibits a dichotomy between \emph{clustering} (mono-type equilibrium) and \emph{coexistence} (multi-type equilibrium). This dichotomy is determined by a clustering criterion that is given in terms of the dual of the system, which consists of a system of \emph{interacting} coalescing random walks. In this paper we provide an alternative clustering criterion, given in terms of an auxiliary dual that is simpler than the original dual, and identify {a} range of parameters for which the criterion is met, which we refer to as the \emph{clustering regime}. It turns out that if the sizes of the active populations are non-clumping (i.e., do not take arbitrarily large values in finite regions of the geographic space) and the relative strengths of the seed-banks {(i.e., the ratio of the sizes of the dormant and the active population in each colony) are bounded uniformly over the geographic space}, then clustering prevails if and only if the symmetrised migration kernel is recurrent.

The spatial system is hard to analyse because of the interaction in the original dual and the inhomogeneity of the colony sizes. By comparing the auxiliary dual with a \emph{non-interacting} two-particle system, we are able to control the correlations that are caused by the interactions. The work in \cite{HN01} and the present paper is part of a broader attempt to include dormancy into interacting particle systems.

\medskip\noindent
\emph{Keywords:} 
Moran model, resampling, migration, seed-bank, inhomogeneity, duality, coexistence versus clustering.

\medskip\noindent
\emph{MSC 2010:} 
Primary 
60K35; 
Secondary 
92D25. 

\medskip\noindent 
\emph{Acknowledgements:} 
The research in this paper was supported by the Netherlands Organisation for Scientific Research (NWO) through grant TOP1.17.019. The authors thank Simone Floreani, Cristian Giardin\`a and Frank Redig for discussions on duality, and Evgeny Verbitskiy for discussions on ergodic theory.   

\end{abstract}

\bigskip

\footnoterule
\noindent
\hspace*{0.3cm} {\footnotesize $^{1)}$ 
Mathematisch Instituut, Universiteit Leiden, Niels Bohrweg 1, 2333 CA  Leiden, NL\\
denholla@math.leidenuniv.nl}\\
\hspace*{0.3cm} {\footnotesize $^{2)}$ 
Mathematisch Instituut, Universiteit Leiden, Niels Bohrweg 1, 2333 CA  Leiden, NL\\
s.nandan@math.leidenuniv.nl}

\newpage
\tableofcontents
\newpage
\section{Introduction}


\subsection{Background and targets}
{\emph{Seed-bank}, as the name suggests, refers to a reservoir that can store genetic information of a biological population for future evolutionary purposes. While this terminology naturally applies to the population of plants, in the literature of population genetics it also relates to \emph{dormancy}, which is a biological trait observed in many microbial communities. The trait provides an organism with the ability to enter a reversible state of low metabolic activity for an indefinite period in unfavourable environmental conditions. In the dormant state an organism refrains from reproduction and other phenotypic developments, until it becomes metabolically active again under better ecological circumstances. Recent studies have revealed that the presence of a seed-bank (or dormancy) significantly changes the way in which a population behaves under evolutionary forces such as migration, selection, mutation or recombination (for references see \cite{LHBB}). Consequently, better understanding dormancy from a biological and a mathematical perspective has attracted the attention of researchers in population genetics. Various attempts have been made to include seed-banks in already existing mathematical models describing genetic evolution of populations.}

In \cite{BCEK15} and \cite{BCKW16}, the discrete-time Fisher-Wright model with \emph{seed-bank} was introduced and analysed. Individuals live in a colony, are subject to \emph{resampling} (i.e., they randomly adopt each others type), and move in and out of the seed-bank, where they suspend their resampling. The seed-bank acts as a repository for the genetic information of the population. Individuals that reside inside the seed-bank are called \emph{dormant}, those that reside outside are called \emph{active}. Both the long-time scaling and the genealogy of the population were analysed in the limit as the size of the colony tends to infinity. The continuous-time version of the model, called Moran model, has the same behaviour. For a recent overview on seed-bank models in population genetics we refer the reader to \cite{BK}. 

{The demography of natural populations is in general structured, in the sense that they admit a \emph{carrying capacity} that is usually imposed by the surrounding habitat. This motivated us in \cite{HN01} to extend existing models with a seed-bank in an \emph{inhomogeneous spatial} setting. Our model is formally described via an \emph{interacting particle system}, and a key characteristic is that no large-colony-size limit is taken. In particular, the constituent populations have preassigned fixed sizes that can be chosen arbitrarily, as long as their spatial growth is not too irregular. Individuals live in multiple colonies, labelled by $\Z^d$, $d\geq 1$, playing the role of a \emph{geographic space} and carry one of the two genetic types: $\heartsuit$ and $\spadesuit$.} Each colony has its own seed-bank, and individuals resample not only from their own colony, but also from other colonies according to a random walk transition kernel, which is referred to as \emph{migration}. The sizes of the active and the dormant population depend on the colony. It was shown that, under mild conditions on the sizes, the system is well-defined, has a unique equilibrium that depends on the initial density of types, and exhibits a dichotomy between \emph{clustering} (mono-type equilibrium) and \emph{coexistence} (multi-type equilibrium). This dichotomy is determined by a \emph{clustering criterion} that is given in terms of the \emph{dual} of the system, which consists of a system of \emph{interacting} coalescing random walks in an inhomogeneous environment. 

{In the present paper, we identify a subdomain of the \emph{clustering regime} (i.e., a range of parameters for which the clustering criterion is met) that is natural and adequate from a biological point of view}. More precisely, we show that if the sizes of the active populations are non-clumping, i.e., do not take arbitrarily large values in finite regions of the geographic space, and the relative strengths of the seed-banks in the different colonies are bounded, then the dichotomy between coexistence and clustering is the classical dichotomy between transience and recurrence of the symmetrised migration kernel, a property that is known to hold for colonies without seed-bank.

Three open problems for the future are: 
\begin{itemize}
	\item[(A)]
	Identify the clustering regime when the relative strengths of the seed-banks in the different colonies are unbounded. In that setting the clustering regime will be different, because it will be driven by a delicate interplay between migration and seed-bank. 
	\item[(B)] 
	In the coexistence regime, identify the \emph{domain of attraction} of the equilibria.
	\item[(C)] 
	In the clustering regime, identify the \emph{growth rate of the mono-type clusters}. 
\end{itemize}
In \cite{HN01} we only showed convergence to equilibrium starting from a family of initial states that are labelled by the initial density of types and that are products of binomial distributions tuned to the inhomogeneity of the relative strengths of the seed-banks.

In \cite{GHO1}, \cite{GHO2}, \cite{GHO3} a \emph{homogeneous} spatial version of the Fisher-Wright model was considered {(i.e., the relative strengths of the seed-banks do not vary across different colonies)}, in the large-colony-size limit. For three different choices of seed-bank, it was shown that the system is well-defined, has a unique equilibrium that depends on the initial density of types, and exhibits a dichotomy between clustering and coexistence. A full description of the clustering regime was obtained. In addition, the finite-systems scheme was established (i.e., how a truncated version of the system behaves on a properly tuned time scale as the truncation level tends to infinity). Moreover, a multi-scale renormalisation analysis was carried out for the case where the colonies are labelled by the hierarchical group. The respective duals for these models are easier, because they are non-interacting and have no inhomogeneity in space. The dual of our model is much harder, which is why our results are much more modest. 


\subsection{Outline}

The paper is organised as follows. In Section~\ref{s.dichotomy} we give a quick definition of the model and state our main theorems about the dichotomy of clustering versus coexistence by identifying the clustering regime for both. In Section~\ref{s.duals} {we provide a different representation (namely, given by a \emph{coordinate process}) of the two-particle dual process associated to} our system introduced in \cite{HN01}, and define {two} auxiliary duals that serve as comparison objects. We relate the coalescence probabilities of the different duals, which leads to a necessary and sufficient criterion for clustering in our system. In Section~\ref{s.coalbd} we prove our main theorems. {In Appendix~\ref{aps-basic-dual} we recall the original representation (given by a \emph{configuration process}) of the two-particle dual from \cite{HN01}, and briefly elaborate on its relation with the alternative representation given in Section~\ref{s.duals}.}


\section{Main theorems}
\label{s.dichotomy}

In Section~\ref{ss.quick} we give a quick definition of the multi-colony system. In Section~\ref{ss.dich} we state our results about the dichotomy of clustering versus coexistence, which requires additional conditions on the sizes of the active and the dormant population.


\subsection{Quick definition of the system} 
\label{ss.quick}

Individuals live in colonies labelled by $\Z^d$, $d \geq 1$, which {play} the role of a \emph{geographic space}. Each colony has an \emph{active} population and a \emph{dormant} population. For $i\in\Z^d,$ we write $(N_i,M_i)\in\N^2$ to denote the \emph{size} of the active, respectively, dormant population at the colony $i$. Each individual carries one of two \emph{types}: $\heartsuit$ and $\spadesuit$. Individuals are subject to: 
\begin{itemize}
	\item[(1)]
	Active individuals in any colony \emph{resample} with active individuals in any colony. 
	\item[(2)]
	Active individuals in any colony \emph{exchange} with dormant individuals in the same colony.
\end{itemize} 
For (1) we assume that each active individual at colony $i$ at rate $a(i,j)$ uniformly draws an active individual at colony $j$ and \emph{adopts its type}. For (2) we assume that each active individual at colony $i$ at rate $\lambda \in (0,\infty)$ uniformly draws a dormant individual at colony $i$ and the two individuals \emph{trade places while keeping their type} (i.e., the active individual becomes dormant and the dormant individual becomes active). Although the exchange rate $\lambda$ could be made to vary across colonies, for the sake of simplicity we choose it to be constant, and we let the \emph{migration kernel} $a(\cdot\,,\cdot)$ be translation invariant and irreducible. Note that dormant individuals do \emph{not} resample.

\begin{assumption}{\bf [Homogeneous migration]}
	\label{assumpt1}
	The migration kernel $a(\cdot\,,\cdot)$ satisfies:
	\begin{itemize}
		\item 
		$a(\cdot\,,\cdot)$ is irreducible in $\Z^d$.
		\item 
		$a(i,j) = a(0,j-i)$ for all $i,j\in{\Z^d}$.
		\item  
		$c:=\displaystyle\sum_{i\in\Z^d\backslash\{0\}} a(0,i) < \infty$ and $a(0,0)=\frac{1}{2}$.
	\end{itemize}
\end{assumption}

\noindent
The second {part of the} assumption ensures that the way genetic information moves between colonies is homogeneous in space. The third {part of the} assumption ensures that the total rate of resampling {of a single individual} is finite and that resampling is possible also at the same colony. 

{
	\begin{remark}
		\rm A detailed description of the multi-colony system can be found in \cite[Section 3.2]{HN01}. In what follows, the geographic space, which here is chosen to be $\Z^d$, can be any countable Abelian group. Moreover, the choice of $a(0,0)=\tfrac{1}{2}$ in Assumption~\ref{assumpt1} has been made only to make our model fit with the classical single-colony Moran model. The value of $a(0,0)$ represents the rate at which individuals resample from their own colony and in principle can be set to any positive real number.\hfill$\Box$
	\end{remark}
} 

Furthermore, in order to avoid trivial statement we assume the following:

\begin{assumption}{\bf [Non-trivial colony sizes]}
	\label{assump:non-trivial colony-sizes}
	In each colony, both the active and the dormant population consist of at least two individuals, i.e., $N_i\geq 2$ and $M_i\geq 2$ for all $i\in\Z^d$.
\end{assumption}

\noindent
For colony sizes where Assumption~\ref{assump:non-trivial colony-sizes} fails, all the results stated below can be obtained with minor technical modifications.

At each colony $i$ we register the pair $(X_i(t),Y_i(t))$, representing the number of active, respectively, dormant individuals of type $\heartsuit$ at time $t$ at colony $i$. The resulting Markov process is denoted by
\begin{equation}
	\label{process}
	Z:=(Z(t))_{t \geq 0}, \qquad Z(t) = ((X_i(t),Y_i(t))_{i \in \Z^d},
\end{equation}
and lives on the state space
\begin{equation}
	\mathcal{X} = \prod_{i\in\Z^d} [N_i]\times[M_i],
\end{equation}
where $[n] := \{0,1,\ldots,n\}$, $n \in \N$. In \cite{HN01}, it was shown that under mild assumptions on the model parameters, the Markov process in \eqref{process} is well defined and has a \emph{dual} $(Z^*(t))_{t \geq 0}$ where the process
\begin{equation}
	Z^*:=(Z^*(t))_{t\geq 0},\quad Z^*(t):= (n_i(t),m_i(t))_{i\in\Z^d},
\end{equation}
lives on the state space
\begin{equation}
	\X^* := \Big\{(n_i,m_i)_{i\in\Z^d} \in \mathcal{X}\colon\,\sum_{i\in\Z^d}(n_i+m_i)<\infty\Big\}.
\end{equation}
The dual process $Z^*$ consists of finite collections of particles that switch between an \emph{active} state and a \emph{dormant} state, and perform \emph{interacting coalescing} random walks while in the active state, with rates that are controlled by the model parameters. 

We recall below the results from \cite{HN01} on the well-posedness of the process $Z$ and the duality relation between $Z$ and $Z^*$. {Interested readers can find the precise duality relation between $Z$ and $Z^*$ in \cite[Theorem 3.10]{HN01}}.

\begin{theorem}{\bf [Well-posedness and duality]{\rm \cite[Theorem 2.2 and Corollary 3.11]{HN01}}}
	\label{T.WPD}
	Suppose that Assumption~{\rm \ref{assumpt1}} is in force. Then the martingale problem associated with \eqref{process} is well-posed under either of the two following conditions:
	\begin{itemize}
		\item[\rm{(a)}]
		$\lim_{\|i\| \to \infty} \|i\|^{-1} \log N_i = 0$ and $\sum_{i\in\Z^d} \e^{\delta \|i\|} a(0,i)<\infty$ for some $\delta>0$.
		\item[\rm{(b)}]
		$\sup_{i\in\Z^d\backslash\{0\}} \|i\|^{-\gamma} N_i < \infty$ and $\sum_{i\in\Z^d} \|i\|^{d+\gamma+\delta} a(0,i)<\infty$ for {some} $\gamma>0$ and some $\delta>0$.
	\end{itemize}
	Furthermore, the Markov process $(Z(t))_{t\geq 0}$ has a factorial moment dual $(Z^*(t))_{t\geq 0}$, living on the state space $\mathcal{X}^*\subset \mathcal{X}$ and consisting of all configurations with finite mass.
\end{theorem}

\noindent
In view of the above, from here onwards, we implicitly assume that the model parameters $(N_i)_{i\in\Z^d}$ and $a(\cdot\,,\cdot)$ are such that one of the two conditions (a) and (b) is satisfied.

{
	\begin{remark}
		\label{remark:moment-regularity}
		{\rm Due to conditions (a) and (b) in Theorem~\ref{T.WPD}, the migration kernel $a(\cdot\,,\cdot)$ always admits a $(d+\delta)$-moment for some $\delta>0$.}\hfill$\Box$
	\end{remark}
}\noindent
We write $\hat{a}(\cdot\,,\cdot)$ to denote the \emph{symmetrised migration kernel} defined by
\begin{equation}
	\hat{a}(i,j):=\tfrac{1}{2}[a(i,j)+a(j,i)], \qquad i,j\in\Z^d,
\end{equation}
and write $a_n(\cdot\,,\cdot)$ to denote the $n$-step transition probability kernel of the embedded chain associated to the continuous-time random walk on $\Z^d$ with rates $a(\cdot\,,\cdot)$. Furthermore, we denote by $\hat{a}_t(\cdot\,,\cdot)$ (respectively, $a_t(\cdot\,,\cdot)$), the time-$t$ transition probability kernel of the continuous-time random walk with migration rates $\hat{a}(\cdot\,,\cdot)$ (respectively, $a(\cdot\,,\cdot)$), and put
\begin{equation}
	K_i:=\frac{N_i}{M_i}, \qquad i\in\Z^d, 
\end{equation}
for the \emph{ratios} of the sizes of the active and the dormant population in each colony. Note that $K_i^{-1}$ quantifies the \emph{relative strength} of the seed-bank at colony $i\in\Z^d$. 

Let $\mathcal{P}$ be the set of probability distributions on $\mathcal{X}$ defined by
\begin{equation} 
	\mathcal{P} := \big\lbrace\mathcal{P}_\theta\colon\,\theta \in [0,1]\big\rbrace,
	\qquad \mathcal{P}_\theta := \theta \prod_{i\in\Z^d}\delta_{(0,0)} + (1-\theta) \prod_{i\in\Z^d}\delta_{(N_i,M_i)}.
\end{equation}
We say that \eqref{process} exhibits \emph{clustering} if the limiting distribution of $Z(t)$ (given that it exists) falls in $\mathcal{P}$. Otherwise we say that it exhibits \emph{coexistence}. In the next section we state the \emph{clustering criterion} from \cite{HN01}, given in terms of the original dual process $Z^*$, and provide an alternative equivalent criterion in terms of a simpler two-particle process that is absorbing.


\subsection{Clustering versus coexistence}
\label{ss.dich}

In \cite{HN01}, it was shown that the system admits a mono-type equilibrium (clustering) if and only if the following criterion is met:

\begin{theorem}{\bf [Clustering condition]}{\rm{\cite[Theorem 3.17]{HN01}}}\label{T.Dichotomy}
	The system clusters if and only if in the dual process $Z^*$ two particles, starting from any locations in $\Z^d$ and any states (active or dormant), coalesce with probability $1$.
\end{theorem}

Before we state our alternative criterion for clustering, we introduce an auxiliary two-particle dual process. In Proposition~\ref{prop:stability of two-particle dual}, we will show the well-posedness of this process. Recall that $\lambda$ is the exchange rate between active and dormant individuals in each colony.

\begin{definition}{\bf [Auxiliary two-particle system]}
	\label{definition: interacting RW2}
	{\rm The two-particle process 
		$\hat\xi := (\hat\xi(t))_{t \geq 0}$ is a continuous-time Markov chain on the state space 
		\begin{equation}
			\mathcal{S}:=(G\times G)\cup\{\circledast\},\quad G:=\Z^d\times\{0,1\}
		\end{equation} 
		with transition rates
		\begin{equation}
			\label{eqn:aux-two-particle-transition-rate}
			\begin{aligned}
				&[(i,\alpha),(j,\beta)] \to\\
				&\begin{cases}
					\displaystyle
					\circledast, &\text{ at rate } \displaystyle
					{2a(i,i)}\tfrac{\alpha\beta}{N_i}\delta_{i,j},\\				
					[(i,1-\alpha),(j,\beta)], &\text{ at rate } 
					\lambda[\alpha+(1-\alpha)K_i]-\tfrac{\lambda}{M_i}\delta_{i,j}(1-\delta_{\alpha,\beta}),\\
					[(i,\alpha),(j,1-\beta)], &\text{ at rate } 
					\lambda[\beta+(1-\beta)K_j]-\tfrac{\lambda}{M_j}\delta_{i,j}(1-\delta_{\alpha,\beta}),\\
					[(k,\alpha),(j,\beta)], &\text{ at rate } 
					\alpha\, a(i,k) \quad \,\text{ for } k\neq i\in\Z^d,\\
					[(i,\alpha),(k,\beta)], &\text{ at rate } 
					\beta\, a(j,k) \quad \,\text{ for } k\neq j\in\Z^d,
				\end{cases}
			\end{aligned}
		\end{equation}
		where $[(i,\alpha),(j,\beta)]\in G\times G$ and  $\delta_{\cdot,\cdot}$ denotes the Kronecker delta-function.} \hfill $\Box$
\end{definition}

\noindent
Here, $\hat{\xi}(t)=[(i,\alpha),(j,\beta)]$ {captures} the location ($i,j \in \Z^d$) and the state ($\alpha,\beta\in\{0,1\}$) of the two particles at time $t$, where $0$ stands for dormant and $1$ stands for active, respectively. Note that $\circledast$ is an absorbing state for the process $\hat\xi$, which is absorbed at a location-dependent rate only when the two particles are on top of each other and in the active state. {We will see in Section~\ref{ss.basicduals} that this is different from what happens in the two-particle system obtained from the original dual}. The process $\hat \xi$ is much simpler than the original two-particle system, because the particles do not interact unless they are on top of each other with opposite states. {Indeed, note that in the second and third line of \eqref{eqn:aux-two-particle-transition-rate} the second term represents a repulsive interaction between the two particles that is non-zero only when $i=j$ and $\alpha\neq \beta$}. From here onwards, we write $\hat{\P}^\eta$ to denote the law of the process $\hat{\xi}$ started from $\eta\in\mathcal{S}$, and $\hat{\E}^\eta$ to denote expectation w.r.t.\ $\hat{\P}^\eta$.

\begin{remark}
	\label{remark: irreducible hat system}
	{\rm Note that, by virtue of {Assumptions~\ref{assumpt1} and \ref{assump:non-trivial colony-sizes}}, all states in $\mathcal{S}$ are accessible by $\hat\xi$.}\hfill$\Box$
\end{remark}

\begin{theorem}{\bf[Clustering criterion]}
	\label{T.Alt-Dichotomy}
	The system clusters if the process $\hat{\xi}$ starting from an arbitrary configuration in $G\times G$ is absorbed with probability $1$. Furthermore, if the sizes of the active populations are non-clumping, i.e.,
	\begin{equation}
		\label{eqn:non-clumping criterion}
		\inf_{i\in\Z^d} \sum_{\|j-i\|\leq R} \tfrac{1}{N_j}>0 \text{ for some } R<\infty, 
	\end{equation} 
	then the converse is true as well.
\end{theorem}

\begin{remark}
	{\rm The condition in \eqref{eqn:non-clumping criterion} is equivalent to requiring that, for some constant $C<\infty$ and all $i\in\Z^d$, there exists a $j$ with $\|j-i\| \leq R$ such that $N_j\leq C$. This requirement can be further relaxed to
		\begin{equation}
			\inf_{i\in\Z^d} \sum_{j\in\Z^d} \frac{1}{N_j} \sum_{n\in\N} m^{2n}\,a_n(i,j)^2 > 0,
		\end{equation}
		where $m:=\frac{c}{2(c+\lambda)+1}$. {Although \eqref{eqn:non-clumping criterion} arises in our context as a technical requirement, it has an interesting connection with the notion introduced in \cite{NK2002, SKLN05} of \emph{coalescent effective population size} (CES) in a subdivided population. Roughly, $N \in \N$ is said to be the CES of a subdivided population when, after measuring time in units of $N$ generations and taking the large-colony-size-limit, the associated genealogy gives rise to Kingman's coalescent (or a similar object). When migration is controlled by a transition matrix, the CES is often proportional to the harmonic mean of the constituent population sizes (see e.g., \cite{W2001}, and also \cite[Section 4.4]{Durrett08}). The non-clumping criterion in \eqref{eqn:non-clumping criterion} essentially says that if $H(i,R)$ is the harmonic mean of the active population sizes of the colonies within the $R$-neighbourhood of colony $i$, i.e.,
			\begin{equation}
				H(i,R) := \frac{|\{j\in\Z^d\,:\,\|j-i\|\leq R\}|}{\displaystyle\sum_{\|j-i\|\leq R} N_j^{-1}}, \qquad i\in\Z^d,
			\end{equation} 
			then $\sup_{i\in\Z^d} H(i,R)<\infty$ for some $R<\infty$. We believe that this connection of the non-clumping criterion to the CES is not accidental, and merits further investigation.
			\hfill$\Box$}}
\end{remark}
To verify when the above clustering criterion is satisfied, we need to impose the following regularity condition on the migration kernel.

{
	\begin{assumption}{\bf [Regularly varying migration kernel]}
		\label{assumpt2} 
		Assume that $t\mapsto\hat{a}_{t}(0,0)$ is regularly varying at infinity, i.e., $\lim_{t\to\infty} \frac{\hat{a}_{pt}(0,0)}{\hat{a}_t(0,0)}=p^{-\sigma}$ for all $p \in (0,\infty)$ and some $\sigma \in [0,\infty)$, where $-\sigma$ is the index of the regular variation and $\hat{a}_t(\cdot,\cdot)$ is the time-$t$ symmetrised migration kernel.
	\end{assumption}
}

{
	\begin{remark}
		{\rm	Note that all genuinely $d$-dimensional continuous-time random walks satisfying the LCLT (see e.g., \cite[Chapter 2]{Lawler2010}) have a probability transition kernel with a regularly varying tail of index $-\tfrac{d}{2}$.\hfill$\Box$}
	\end{remark}
}
When the relative strengths of the seed-banks are uniformly bounded, clustering is equivalent to the symmetrised migration kernel being recurrent, a setting that is classical. The following theorem provides a slightly weaker result. 

\begin{theorem}{\bf [Clustering regime]}
	\label{T.Coal_recc}
	Suppose that Assumption~\ref{assumpt2} is in force. Assume that the active population sizes are non-clumping, i.e., \eqref{eqn:non-clumping criterion} is satisfied, and the relative strengths of the seed-banks are uniformly bounded, i.e.,
	\begin{equation}
		\label{as: bounded seed-bank}
		\sup_{i\in\Z^d} K_i^{-1} < \infty.
	\end{equation}
	If the system clusters, then it is necessary that the symmetrised kernel $\hat{a}(\cdot\,,\cdot)$ is recurrent. Furthermore, if the migration kernel $a(\cdot\,,\cdot)$ is symmetric, then the converse holds as well.
\end{theorem}

{It was shown in \cite{GHO1} that the above dichotomy is true when the seed-banks are homogeneous (i.e., $(N_i,M_i)=(N,M)$ for all $i\in\Z^d$) and the large-colony-size limit is taken (i.e., $N,M \to \infty$ such that $N/M \to K \in (0,\infty)$). In that case, the dual process is an independent particle system with coalescence and without inhomogeneity, for which the proof is much simpler. The result stated above extends the dichotomy to the inhomogeneous setting. It essentially says that if the inhomogeneities caused by the seed-banks are spatially uniform (reflected by \eqref{as: bounded seed-bank}), then the dichotomy remains unchanged. The condition in \eqref{as: bounded seed-bank} allows us to compare the auxiliary dual process $\hat{\xi}$ with a non-interacting two-particle process $\xi^*$ living on the state space $\mathcal{S}$ that we introduce in Section~\ref{ss.basicduals} (see Section~\ref{ss.indep-sys} for more details). As we will see later, $\circledast$ is an absorbing state for $\xi^*$ too and, under the conditions given in Theorem~\ref{T.Coal_recc}, it turns out that $\hat{\xi}$ is absorbed with probability 1 if and only if $\xi^*$ is. In $\xi^*$ the two particles evolve independently until absorption. A single particle migrates in the active state at rates $a(\cdot\,,\cdot)$, becomes dormant from the active state at rate $\lambda$, and becomes active from the dormant state at rate $\lambda K_i$ when it is at location $i\in\Z^d$. When the condition \eqref{as: bounded seed-bank} is met, the average time spent in the dormant state by the particles in the various locations are of the same order, and hence the distance between the two particles is effectively controlled by the symmetrised kernel $\hat{a}(\cdot\,,\cdot)$. In particular, the recurrence of $\hat{a}(\cdot\,,\cdot)$ forces the two particles to meet each other infinitely often with probability 1. As a result, $\xi^*$ is eventually absorbed in $\circledast$. We  exploit these facts along with the alternative clustering criterion to prove Theorem~\ref{T.Coal_recc}. We expect the symmetry assumption to be redundant for the converse statement, but are unable to remove it for technical reasons. The following result is an immediate corollary.}

	\begin{corollary}{\bf [Dimensional dichotomy]}
		\label{C.dichotomy}
		Assume that all the conditions in Theorem~\ref{T.Coal_recc} are in force. Then the following hold:
		\begin{enumerate}[{\rm (a)}]
			\item 
			Coexistence prevails when $d>2$.
			\item 
			Clustering prevails when $d\leq 2$ and $a(\cdot\,,\cdot)$ is symmetric. 
		\end{enumerate}
	\end{corollary}	


\section{Dual processes: comparison between different systems}
\label{s.duals}

{In Section~\ref{ss.basicduals} we give a brief description of the dual process $Z^*$ of our original system introduced in \cite{HN01}, and define \emph{two auxiliary duals} that serve as comparison objects. The auxiliary duals are simplified versions of the basic dual, started from two particles, where the coalesced state of the two particles is turned into an absorbing state. In Sections~\ref{ss.comp12}--\ref{ss.compindep} we relate the coalescence (absorption) probabilities of the auxiliary duals via a comparison technique that is based on the Lyapunov function approach employed in \cite{CK11}. In Section~\ref{ss.concl} we provide finer conditions on the parameters of our original model under which the results derived in previous sections hold.}  


\subsection{Two-particle dual and auxiliary duals}
\label{ss.basicduals} 

{
	Recall that the dual process $Z^*$ is an interacting particle system describing the evolution of finitely many particles such that (see \cite[Section 3.2]{HN01} for more details)
	\begin{itemize}
		\item particles can be in one of the two states: \emph{active} and \emph{dormant},
		\item particles migrate while in the active state,
		\item a pair of particles in the active state can coalesce (even from different locations) with each other to form a single active particle,
		\item the interaction between the particles is \emph{repulsive} in nature, in the sense that a particle discourages another particle to be at the same location with the same state (active or dormant). To be more precise, the associated transition of a particle happens at a slower rate due to the interaction with the other particles.
	\end{itemize}
	As stated earlier in Theorem~\ref{T.Dichotomy}, the dichotomy between clustering and coexistence is solely determined by the coalescence of two dual particles, and so we only need to analyse the dual process starting from two particles. There are two ways in which we can describe the two-particle dual process, namely, as a \emph{configuration process} that keeps track of the number of active and dormant particles at each location of the geographic space, or as a \emph{coordinate process} that gives only the location and the state (active or dormant) of the particles that are present in the system. In \cite{HN01}, the dual process $Z^*$ was introduced via a configuration process. However, in what follows we describe the two-particle dual originating from the process $Z^*$ as a coordinate process in order to keep computations and notations simple. To make our paper self-contained, in Appendix~\ref{aps-basic-dual} we include a short description of the configuration process associated to the original two-particle dual.}
{
	The transition rates for the two particles in the dual process are as follows:
	\begin{itemize}
		\item \textbf{(Migration)} 
		An active particle at site $i$ migrates to site $j$ at rate $a(i,j)$ if there is no active particle at site $j$, otherwise at rate $a(i,j)(1-\tfrac{1}{N_j})$.
		\item \textbf{(Active to Dormant)} 
		An active particle at site $i$ becomes dormant at site $i$ at rate $\lambda$ if there is no dormant particle at site $i$, otherwise at rate $\lambda(1-\tfrac{1}{M_i})$.
		\item \textbf{(Dormant to Active)}  
		A dormant particle at site $i$ becomes active at site $i$ at rate $\lambda K_i$ if there is no active particle at site $i$, otherwise at rate $\lambda (K_i-\tfrac{1}{M_i})$.
		\item \textbf{(Coalescence)} 
		An active particle at site $i$ coalesces with another active particle at site $j$ at rate $\tfrac{a(i,j)}{N_j}$.
	\end{itemize}
	Since we are interested in the coalescence probability of two dual particles only, it suffices to let the two-particle process absorb into a single state $\circledast$ as soon as coalescence happens and analyse the absorption time of the resulted absorbing process. Note that in the two-particle dual described above, once coalescence has occurred only a single particle remains in the system for the rest of the time. Because of this, $\circledast$ becomes an absorbing state. Furthermore, by virtue of the well-known Dynkin criterion for lumpability, the absorbing process remains a continuous-time Markov chain. Although this can be verified by standard computations, for the convenience of the reader we include a brief proof in Appendix~\ref{aps-basic-dual}. Below we provide a formal definition of the absorbing two-particle process as \emph{interacting RW1}, which is basically a coordinate process living on the state space
	\begin{equation}
		\mathcal{S}:= (G \times G)\cup\{\circledast\},\quad G := \Z^d\times\{0,1\}.
	\end{equation}
}

\begin{definition}{\bf [Interacting RW1]}
	\label{definition: interacting RW1}
	{\rm The interacting RW1 process 
		\begin{equation}
			\xi := (\xi(t))_{t \geq 0}
		\end{equation} 
		is the continuous-time Markov chain on the state space $\mathcal{S}$ with transition rates
		\begin{equation}
			\label{eqn:two-particle-dual-transition-rates}
			\begin{aligned}
				&[(i,\alpha),(j,\beta)] \to\\
				&\begin{cases}
					\displaystyle
					\circledast, &\text{ at rate } \displaystyle
					\alpha\beta(1-\delta_{i,j})\Big[\tfrac{a(i,j)}{N_j}+\tfrac{a(j,i)}{N_i}\Big]+{2a(i,i)}\tfrac{\alpha\beta}{N_i}\delta_{i,j},\\
					[(i,1-\alpha),(j,\beta)], &\text{ at rate } 
					\lambda[\alpha+(1-\alpha)K_i]-\tfrac{\lambda}{M_i}\delta_{i,j}(1-\delta_{\alpha,\beta}),\\
					[(i,\alpha),(j,1-\beta)], &\text{ at rate } 
					\lambda[\beta+(1-\beta)K_j]-\tfrac{\lambda}{M_j}\delta_{i,j}(1-\delta_{\alpha,\beta}),\\
					[(k,\alpha),(j,\beta)], &\text{ at rate } 
					\alpha\, a(i,k)-a(i,k)\tfrac{\alpha\beta}{N_j}\delta_{k,j} \quad \,\text{ for } k\neq i\in\Z^d,\\
					[(i,\alpha),(k,\beta)], &\text{ at rate } 
					\beta\, a(j,k)-a(j,k)\tfrac{\alpha\beta}{N_i}\delta_{k,i} \quad \,\text{ for } k\neq j\in\Z^d,
				\end{cases}
			\end{aligned}
		\end{equation}
		where $\delta_{\cdot,\cdot}$ denotes the Kronecker delta-function.
	} \hfill $\Box$
\end{definition}

\noindent
Here, $\xi(t)=[(i,\alpha),(j,\beta)]$ provides the location ($i,j \in \Z^d$) and the state ($\alpha,\beta\in\{0,1\}$) of the two particles in the process at time $t$, where $0$ stands for dormant and $1$ stands for active, respectively.

{	
	\begin{remark}
		\label{remark:equivalence-of-absorption}
		\rm Note that the coalescence time of the original two-particle dual process becomes the absorption time of $\xi$, and thus the original clustering criterion stated in Theorem~\ref{T.Dichotomy} is equivalent to asking whether or not $\xi$ is absorbed in $\circledast$ with probability 1. However, the negative second terms in the last two transition rates of $\xi$ (see \eqref{eqn:two-particle-dual-transition-rates}) imply that the two particles interact repulsively with each other even when they migrate in the active state. As a consequence, the effective migration kernel of a single particle becomes inhomogeneous in space, and so $\xi$ is much harder to analyse than the auxiliary two-particle dual $\hat{\xi}$ defined in Definition~\ref{definition: interacting RW2}. Another key difference between $\xi$ and $\hat{\xi}$ is that $\hat{\xi}$ has a positive rate of absorption only when both particles are on the same location in the active state. Although it may seem natural that $\xi$ has a higher chance of absorption than $\hat{\xi}$, we will show later via a comparison argument that, under the non-clumping criterion (see \eqref{eqn:non-clumping criterion}) on $(N_i)_{i\in\Z^d}$, if one process enters the absorbing state $\circledast$ with probability 1, then the other process does so too. This ultimately provides us with the alternative criterion for clustering in Theorem~\ref{T.Alt-Dichotomy}.\hfill$\Box$
	\end{remark}
}
From now onwards, we write $\P^\eta$ to denote the law of the process $\xi$ started from $\eta\in\mathcal{S}$, and $\E^\eta$ to denote expectation w.r.t.\ $\P^\eta$.

\begin{remark}
	\label{remark: lives-on-valid-state}
	{\rm Note that, by virtue of {Assumption~\ref{assumpt1} and Assumption~\ref{assump:non-trivial colony-sizes}}, all states in $\mathcal{S}$ are accessible by $\xi$.}\hfill$\Box$
\end{remark}

In addition to the auxiliary two-particle process $\hat{\xi}$ defined in Definition~\ref{definition: interacting RW2}, and  and the \emph{interacting RW1} process $\xi$ defined above, we introduce one more two-particle system, called \emph{independent RW}, on the same state space $\mathcal{S}$. This will also serve as an intermediate comparison object. 

\begin{definition}{\bf [Independent RW]}
	\label{definition: independent RW}
	{\rm The independent RW process 
		\begin{equation}
			\xi^* := (\xi^*(t))_{t \geq 0}
		\end{equation} 
		is the continuous-time Markov chain on the state space $\mathcal{S}$ with transition rates
		\begin{equation}
			\begin{aligned}
				&[(i,\alpha),(j,\beta)] \to\\
				&\begin{cases}
					\displaystyle
					\circledast, &\text{ at rate } \displaystyle
					{2a(i,i)}\tfrac{\alpha\beta}{N_i}\delta_{i,j},\\
					[(i,1-\alpha),(j,\beta)], &\text{ at rate } 
					\lambda[\alpha+(1-\alpha)K_i],\\
					[(i,\alpha),(j,1-\beta)], &\text{ at rate } 
					\lambda[\beta+(1-\beta)K_j],\\
					[(k,\alpha),(j,\beta)], &\text{ at rate } 
					\alpha\, a(i,k) \quad \,\text{ for } k\neq i\in\Z^d,\\
					[(i,\alpha),(k,\beta)], &\text{ at rate } 
					\beta\, a(j,k) \quad \,\text{ for } k\neq j\in\Z^d.
				\end{cases}
			\end{aligned}
		\end{equation}
	} \hfill $\Box$
\end{definition}

\noindent
{In Section~\ref{s.coalbd} we delve deeper into the independent RW process $\xi^*$, in order to utilize the comparison results derived in the next two sections and determine the clustering regime.} We write $\P^*_\eta$ to denote the law of $\xi^*$ started from $\eta\in\mathcal{S}$, and $\E^*_\eta$ to denote expectation w.r.t.\ $\P^*_\eta$.

{In the following proposition, we establish the well-posedness of $\hat{\xi}$ (see Definition~\ref{definition: interacting RW2}), and of $\xi$, $\xi^*$ defined above.}
\begin{proposition}{\bf[Stability]}
	\label{prop:stability of two-particle dual}
	All three processes $\xi$, $\hat{\xi}$, $\xi^*$ are non-explosive continuous-time Markov chains on the countable state space $\mathcal{S}$.
\end{proposition}

\begin{proof}
	We prove this claim by using the Foster-Lyapunov criterion (see \cite{MT93}). Let $B_0 := \{\circledast\}$, and for $n\in\N$ define  $B_n:= \{[(i,\alpha),(j,\beta)]\in\mathcal{S}\colon\,\max\{\|i\|,\|j\|\}< n\}\cup B_0$. Define
	\begin{equation}
		\label{eqn: lyapunov function}
		\begin{aligned}
			V(\eta):=\begin{cases}
				\|i\|+\|j\|, & \text{ if }\eta=[(i,\alpha),(j,\beta)],\\
				0, &\text{ otherwise,}
			\end{cases}
			\qquad
			\eta\in\mathcal{S}.
		\end{aligned}
	\end{equation}
	Furthermore, let $Q,\hat{Q},Q^*$ be the infinitesimal generators of the processes $\xi,\hat{\xi},\xi^*$, respectively. Note that, for $\eta=[(i,\alpha),(j,\beta)]\in\mathcal{S}$,
	\begin{equation}
		\label{eqn: lyapunov-criterion}
		\begin{aligned}
			QV(\eta) = \alpha&\sum_{k\neq i}a(i,k)(\|k\|-\|i\|)+\beta\sum_{k\neq j}a(j,k)(\|k\|-\|j\|)\\
			&\qquad\qquad  -2\alpha\beta(1-\delta_{i,j})\Big[\tfrac{a(i,j)}{N_j}\|i\|+\tfrac{a(j,i)}{N_i}\|j\|\Big]-\tfrac{\alpha\beta}{N_i}\delta_{i,j}\\
			&\leq (\alpha+\beta)\mu_1+2\alpha\beta(1-\delta_{i,j})\Big[\tfrac{a(i,j)}{N_j}\|i\|+\tfrac{a(j,i)}{N_i}\|j\|\Big]\\
			&\leq 2V(\eta)+(\alpha+\beta)\mu_1\\
			&\leq 2V(\eta)+2\mu_1\quad\text{(since $\alpha+\beta\leq 2$)},
		\end{aligned}
	\end{equation}
	where $\mu_1:= \sum_{i\in\Z^d/\{0\}}\|i\|\,a(0,i)$. Let $V^\prime\colon\,\mathcal{S}\to[0,\infty)$ be the function defined by $\eta\mapsto V(\eta)+\mu_1$. Note that $B_n\uparrow\mathcal{S}$ as $n\to\infty$ and $\inf_{\eta\in B_n^c}V^\prime(\eta)\geq n$. Thus, $\inf_{\eta\in B_n^c}V^\prime(\eta)\uparrow\infty$ as $n\to\infty$ and, by \eqref{eqn: lyapunov-criterion}, $QV^\prime(\eta)\leq 2\,V^\prime(\eta)$. Hence the Foster-Lyapunov criterion is satisfied by the generator $Q$, and so $\xi$ is non-explosive. Similar arguments show that $\hat{\xi}$ and $\xi^*$ are non-explosive as well.
\end{proof}


\subsection{Comparison between interacting duals}

{In this section we show, via comparison between the infinitesimal generators of the two-particle dual $\xi$ and the auxiliary two-particle dual $\hat{\xi}$ introduced in Definition~\ref{definition: interacting RW2}, that the two processes have in fact very similar behaviour when it comes to long-run survivability. This is not surprising given that there are only slight differences in the migration and absorption mechanism (cf.\ the first and the last transition rates in Definition~\ref{definition: interacting RW2} and Definition~\ref{definition: interacting RW1}) of the active particles present in the two processes.}
\label{ss.comp12}

\begin{proposition}{\bf [Stochastic domination]}
	\label{proposition: stochastic domination}
	Let $f\colon\,\mathcal{S}\to\R$ be bounded and such that $f(\eta)\leq f(\circledast)$ for all $\eta \in\mathcal{S}$. Let $(\xi(t))_{t\geq 0}$ and $(\hat\xi(t))_{t\geq 0}$ be the interacting RW1 and the auxiliary two-particle system defined in Definition~\ref{definition: interacting RW1} and Definition~\ref{definition: interacting RW2}, respectively. Then, for any $\eta\in\mathcal{S}$ and $t\geq 0$, $\E^\eta[f(\xi(t))]\geq \hat{\E}^\eta[f(\hat{\xi}(t))]$.
\end{proposition}

\begin{proof}
	Let $Q$ and $\hat{Q}$ be the generators of the processes $\xi,\hat{\xi}$, respectively. Since $\xi$ and $\hat\xi$ are non-explosive continuous-time Markov processs on a countable state space, $Q$ and $\hat{Q}$ generate unique Markov semigroups $(S_t)_{t\geq 0}$ and $(\hat{S}_t)_{t\geq 0}$, respectively, given by
	\begin{equation}
		(S_tg)(\eta)=\E^{\eta}[g(\xi(t))],\quad (\hat{S}_tg)(\eta) = \hat{\E}^{\eta}[g(\hat{\xi}(t))], \qquad t\geq 0,
	\end{equation}
	where $g\colon\,\S\to\R$ is bounded and $\eta\in\S$. Since $f$ is bounded, we can apply the variation of constants formula for semigroups, to obtain
	\begin{equation}
		\label{eqn: variation-of-constants}
		\begin{aligned}
			(S_tf)(\eta)-(\hat{S}_tf)(\eta) = \int_{0}^{t}(S_{t-s}(Q-\hat{Q})\hat{S}_sf)(\eta)\,\d s.
		\end{aligned}
	\end{equation}
	{The actions of $Q$ and $\hat{Q}$ on a bounded function $g\colon\,\S\to\R$ are given by
		\begin{equation}
			\label{eq:generator-of-xi}
			\begin{aligned}
				Qg(\eta) = \alpha&\sum_{k\neq i}a(i,k)\big[1-\tfrac{\beta}{N_j}\delta_{k,j}\big]\{g([(k,\alpha),(j,\beta)])-g([(i,\alpha),(j,\beta)])\}\\
				&+\beta\sum_{k\neq j}a(j,k)\big[1-\tfrac{\alpha}{N_i}\delta_{k,i}\big]\{g([(i,\alpha),(k,\beta)])-g([(i,\alpha),(j,\beta)])\}\\
				&+\big[\lambda(\alpha+(1-\alpha)K_i)-\tfrac{\lambda}{M_i}\delta_{i,j}(1-\delta_{\alpha,\beta})\big]\{g([(i,1-\alpha),(j,\beta)])-g([(i,\alpha),(j,\beta)])\}\\
				&+\big[\lambda(\beta+(1-\beta)K_j)-\tfrac{\lambda}{M_j}\delta_{i,j}(1-\delta_{\alpha,\beta})\big]\{g([(i,\alpha),(j,1-\beta)])-g([(i,\alpha),(j,\beta)])\}\\
				&+\big[\alpha\beta(1-\delta_{i,j})\big(\tfrac{a(i,j)}{N_j}+\tfrac{a(j,i)}{N_i}\big)+2a(i,i)\tfrac{\alpha\beta}{N_i}\delta_{i,j}\big]\{g(\circledast)-g([(i,\alpha),(j,\beta)])\}
			\end{aligned}
		\end{equation}
		and
		\begin{equation}
			\label{eq:generator-of-xi-hat}
			\begin{aligned}
				\hat{Q}g(\eta) = \alpha&\sum_{k\neq i}a(i,k)\{g([(k,\alpha),(j,\beta)])-g([(i,\alpha),(j,\beta)])\}\\
				&+\beta\sum_{k\neq j}a(j,k)\{g([(i,\alpha),(k,\beta)])-g([(i,\alpha),(j,\beta)])\}\\
				&+\big[\lambda(\alpha+(1-\alpha)K_i)-\tfrac{\lambda}{M_i}\delta_{i,j}(1-\delta_{\alpha,\beta})\big]\{g([(i,1-\alpha),(j,\beta)])-g([(i,\alpha),(j,\beta)])\}\\
				&+\big[\lambda(\beta+(1-\beta)K_j)-\tfrac{\lambda}{M_j}\delta_{i,j}(1-\delta_{\alpha,\beta})\big]\{g([(i,\alpha),(j,1-\beta)])-g([(i,\alpha),(j,\beta)])\}\\
				&+\big[2a(i,i)\tfrac{\alpha\beta}{N_i}\delta_{i,j}\big]\{g(\circledast)-g([(i,\alpha),(j,\beta)])\},
			\end{aligned}
		\end{equation}
		where $\eta=[(i,\alpha),(j,\beta)]\in\mathcal{S}$. Thus,
		\begin{equation}
			\label{eqn: explicit-generator-difference}
			\begin{aligned}
				((Q-\hat{Q})g)(\eta) 
				&= -\alpha\beta\sum_{k\neq i}\tfrac{a(i,k)}{N_j}\delta_{k,j}\{g([(k,\alpha),(j,\beta)])-g([(i,\alpha),(j,\beta)])\}\\
				&\qquad-\alpha\beta\sum_{k\neq j}\tfrac{a(j,k)}{N_i}\delta_{k,i}\{g([(i,\alpha),(k,\beta)])-g([(i,\alpha),(j,\beta)])\}\\
				&\qquad+\big[\alpha\beta(1-\delta_{i,j})\big(\tfrac{a(i,j)}{N_j}+\tfrac{a(j,i)}{N_i}\big)\big]\{g(\circledast)-g([(i,\alpha),(j,\beta)])\}\\
				&=-\alpha\beta(1-\delta_{i,j})\big[\tfrac{a(i,j)}{N_j}g([(j,\alpha),(j,\beta)])+\tfrac{a(j,i)}{N_i}g([(i,\alpha),(i,\beta)])\big]\\
				&\qquad+\alpha\beta(1-\delta_{i,j})\big[\tfrac{a(i,j)}{N_j}+\tfrac{a(j,i)}{N_i}\big]g(\circledast),
			\end{aligned}
		\end{equation}
		and so if $g$ is such that $\sup_{\eta\in\S} g(\eta)=g(\circledast)$, then}
	\begin{equation}
		\label{eqn: positivity}
		\begin{aligned}
			((Q-\hat{Q})g)(\eta) &= 
			\begin{cases}
				\alpha\beta(1-\delta_{i,j})\tfrac{a(i,j)}{N_j}\Big[g(\circledast)-g([(j,\alpha),(j,\beta)])\Big]\\
				\quad + \alpha\beta(1-\delta_{i,j})\tfrac{a(j,i)}{N_i}\Big[g(\circledast)-g([(i,\alpha),(i,\beta)])\Big], 
				&\text{ if }\eta=[(i,\alpha),(j,\beta)]\neq\circledast,\\
				0,&\text{ otherwise,}
			\end{cases}\\
			&\geq 0.
		\end{aligned}
	\end{equation}
	Note that the semigroup $(\hat{S}_t)_{t\geq 0}$ also has the property $\sup_{\eta\in\S}(\hat{S}_sf)(\eta) = f(\circledast) = (\hat{S}_sf)(\circledast)$ for any $s\geq 0$, since $f\leq f(\circledast)$ and $\circledast$ is absorbing. Thus, combining the above with \eqref{eqn: positivity}, we get that $(Q-\hat{Q})\hat{S}_sf$ is a non-negative function for any $s\geq 0$. Therefore the right-hand side of \eqref{eqn: variation-of-constants} is non-negative as well, which proves the desired result.
\end{proof}

\begin{corollary}{\bf [Stochastic ordering of absorption times]}
	\label{corollary: nu > nu_hat}
	Let $\tau$ and $\hat{\tau}$ denote the absorption time of the processes $\xi$ and $\hat\xi$, respectively. Then, for any $\eta\in\mathcal{S}$ and $t>0$,
	\begin{equation}
		\P^\eta(\tau\leq t) \geq \hat\P^\eta(\hat\tau\leq t).
	\end{equation}
\end{corollary}

\begin{proof}
	This follows by applying Proposition~\ref{proposition: stochastic domination} to the function $f = \mathbf{1}_{\{\circledast\}}$ and using that $\circledast$ is absorbing for both $\xi$ and $\hat{\xi}$.
\end{proof}

{The above result tells that the two particles in the process $\xi$ have a higher chance of absorption than in the auxiliary process $\hat{\xi}$. This fits with intuition: two active particles in $\xi$ can coalesce even when sitting at different locations. In the next result we show that two particles in $\hat{\xi}$ have a higher probability of being on top of each other in the active state or being absorbed than in $\xi$. This is essentially due to the extra repulsive interaction that takes place when an active particle in $\xi$ attempts to migrate, which is absent in $\hat{\xi}$.}

{
	\begin{proposition}{\bf[Stochastic ordering of hitting times]}
		\label{prop:hat-process-more-freedom}
		Let $B\subset\S$ be defined as
		\begin{equation}
			\begin{aligned}
				B:=\{[(i,1),(i,1)]\colon\, i\in\Z^d\}\cup\{\circledast\}.
			\end{aligned}
		\end{equation}
		Let $T_B,\hat{T}_B$ denote the first hitting time of the set $B$ for $\xi$ and $\hat{\xi}$, respectively. Then, for all $y\in\mathcal{S}$,
		\begin{equation}
			\label{eqn: hat-process-more-freedom}
			\begin{aligned}
				\hat{\P}^y(\hat{T}_B<\infty) \geq \P^y(T_B<\infty).
			\end{aligned}
		\end{equation}
	\end{proposition}
}

{
	\begin{proof}
		Let $g\colon\,\S\to[0,1]$ and $\hat{g}\colon\,\S\to[0,1]$ be defined as
		\begin{equation}
			\label{eqn: defn-hitting-probability-of-B}
			\begin{aligned}
				g(y): = \P^y(T_B<\infty),\quad\hat{g}(y):=\hat{\P}^y(\hat{T}_B<\infty),\quad y\in\S.
			\end{aligned}
		\end{equation}
		We are required to show that
		\begin{equation}
			\label{eqn: hat-g-domination}
			\begin{aligned}
				\hat{g}(y)\geq g(y) \text{ for any } y\in\S.                   
			\end{aligned}
		\end{equation}
		To that end, let $Q$ and $\hat{Q}$ be the generators of the processes $\xi$ and $\hat{\xi}$, respectively. Applying $Q-\hat{Q}$ to the function $\hat{g}$, we get from \eqref{eqn: positivity} that
		\begin{equation}
			\label{eqn: Qg-zero-everywhere}
			\begin{aligned}
				&(Q\hat{g})(y)-(\hat{Q}\hat{g})(y)\\ 
				&= \begin{cases}
					\alpha\beta(1-\delta_{i,j})\tfrac{a(i,j)}{N_j}\Big\{\hat{g}(\circledast)-\hat{g}([(j,\alpha),(j,\beta)])\Big\}\\
					\quad + \alpha\beta(1-\delta_{i,j})\tfrac{a(j,i)}{N_i}\Big\{\hat{g}(\circledast)-\hat{g}([(i,\alpha),(i,\beta)])\Big\}, 
					&\text{ if }y=[(i,\alpha),(j,\beta)]\neq\circledast,\\
					0, &\text{ otherwise.}
				\end{cases}
			\end{aligned}
		\end{equation}
		By a first-{jump} analysis of $\hat{\xi}$, we have $(\hat{Q}\hat{g})(y) = 0$ for any $y\notin B$ and $\hat{g}\equiv 1$ on $B$. Thus, the right-hand side of \eqref{eqn: Qg-zero-everywhere} is always 0, and so $(Q\hat{g})(y)=(\hat{Q}\hat{g})(y)= 0$ for any $y\notin B$. Let $y\in\S$ be fixed, and let $\xi$ be started from $y$. Since $\hat{g}$ is bounded and $\xi$ is non-explosive, the process $(M_t)_{t\geq 0}$ defined by $M_t := \hat{g}(\xi(t)) - \int_{0}^{t} (Q\hat{g})(\xi(s))\,\d s$ is a martingale under the law $\P^y$ w.r.t.\ the natural filtration associated to the process $\xi$. Hence the stopped process $(M_{t\wedge T_B})_{t\geq 0}$ is also a martingale. Note that, since $Q\hat{g}=0$ outside $B$, we have $\int_{0}^{t\wedge T_B} (Q\hat{g})(\xi(s))\,\d s=0$ for any $t\geq 0$. Hence $M_{t\wedge T_B} = \hat{g}(\xi({t\wedge T_B}))$ for any $t\geq 0$. By the martingale property, for any $t> 0$,
		\begin{equation}
			\begin{aligned}
				\hat{g}(y)=\E^y[\hat{g}(\xi(0))] = \E^y[\hat{g}(\xi({t\wedge T_B}))] 
				\geq \E^y[\hat{g}(\xi({T_B}))\mathbf{1}_{T_B< t}] = \P^y(T_B<t).
			\end{aligned}
		\end{equation}
		Letting $t\to\infty$, we get $\hat{g}(y)\geq \P^y(T_B<\infty)=g(y)$, which proves \eqref{eqn: hat-process-more-freedom}.	
	\end{proof}
}

{With the help of the above proposition, we can compare the probability of absorption for $\xi$ and $\hat{\xi}$. Corollary~\ref{corollary: nu > nu_hat} implied that $\xi$ is more likely to get absorbed at $\circledast$ than $\hat{\xi}$. The following result, however, tells that, under a certain condition, if $\xi$ is absorbed with probability 1, then so is $\hat{\xi}$.}

\begin{theorem}{\bf [Comparison of absorption probabilities]}
	\label{theorem: nu_hat > nu}
	Let $\nu\colon\,\mathcal{S}\to[0,1]$ and $\hat\nu\colon\,\mathcal{S}\to[0,1]$ be defined by
	\begin{equation}
		\nu(\eta) := \P^\eta(\tau<\infty),\quad\hat{\nu}(\eta):=\hat{\P}^\eta(\hat{\tau}<\infty),
	\end{equation}
	i.e., $\nu(\eta)$ and $\hat{\nu}(\eta)$ are the absorption probabilities of the processes $\xi$ and $\hat{\xi}$, respectively, started from $\eta$. Assume that
	\begin{equation}
		\label{as: uniformly-bounded-nu-hat-probability}
		\begin{aligned}
			\inf\{\hat{\nu}([(i,1),(i,1)])\colon\,i\in\Z^d\}>0.
		\end{aligned}
	\end{equation}
	{For all $\eta\in\mathcal{S}$, if $\nu(\eta) = 1$, then $\hat{\nu}(\eta)=1$.}
\end{theorem}

\begin{proof}
	The proof is by contradiction. If $\eta=\circledast$, then the claim is trivial. So assume that $\hat{\nu}(\eta) < 1$ and $\nu(\eta)=1$ for some $\eta\neq\circledast$. Note that, by the strong Markov property,
	\begin{equation}
		\label{eqn: zero-infimum}
		\begin{aligned}
			\inf_{y\in\S}\hat{\nu}(y)=0.
		\end{aligned}
	\end{equation}
	Moreover, since by Remark~\ref{remark: lives-on-valid-state} the process $\xi$ started from $\eta$ can visit any configuration $y\in\S$ in finite time with positive probability, we have 
	\begin{equation}
		\label{eqn: constant-nu}
		\begin{aligned}
			\nu(y) = 1 \qquad \forall\,y\in\S.
		\end{aligned}
	\end{equation}
	We will show that \eqref{eqn: zero-infimum} and \eqref{eqn: constant-nu} are contradictory. 
	
	{For $y\in\mathcal{S}$, set
		\begin{equation}
			\begin{aligned}
				g(y): = \P^y(T_B<\infty),\quad\hat{g}(y):=\hat{\P}^y(\hat{T}_B<\infty),\quad y\in\S,
			\end{aligned}
		\end{equation}
		where $T_B,\hat{T}_B$ are the hitting times of the set $B:=\{[(i,1),(i,1)]\,:\,i\in\Z^d\}\cup \{\circledast\}$ for $\xi$ and $\hat{\xi}$, respectively.} Now, since $T_B\leq \tau$ a.s., we have $g(y) \geq \nu(y)$ for any $y\in\S$, and combined with \eqref{eqn: constant-nu} this implies that $g\equiv 1$ on $\S$. So by Proposition~\ref{prop:hat-process-more-freedom}, we have
	\begin{equation}
		\label{eqn: constant-g-hat}
		\begin{aligned}
			\hat{g}(y) = \hat{\P}^y(\hat{T}_B<\infty) = 1\text{ for all } y\in\S,
		\end{aligned}
	\end{equation}
	i.e., the process $\hat{\xi}$ started from any configuration $y\in\S$ enters $B$ with probability $1$. Let $\hat{T}$ be the hitting time of the set $\hat{B}:=B\backslash\{\circledast\}$ for the process $\hat{\xi}$, and let 
	\begin{equation}
		\label{eqn: epsilon-defn}
		\begin{aligned}
			\epsilon:=\inf\{\hat{\nu}(y)\colon\,y\in \hat{B}\}.
		\end{aligned}
	\end{equation}
	By \eqref{as: uniformly-bounded-nu-hat-probability}, we have $\epsilon>0$. Note that $\hat{T}\leq \hat{\tau}$ a.s.\ for the process $\hat{\xi}$, since two particles coalesce only when they are on top of each other and are both active, and so $\hat{T}_B = \hat{T}\wedge\hat{\tau} = \hat{T}$ a.s. Therefore, by \eqref{eqn: constant-g-hat}, $\hat{\P}^y(\hat{T}<\infty)=1$ for any $y\in\S$. Therefore, for $y\in\S$,
	\begin{equation}
		\begin{aligned}
			\hat{\nu}(y) &= \hat{\P}^y(\hat{\tau}<\infty) =\hat{\P}^y(\hat{T}\leq \hat{\tau}<\infty)
			=\sum_{x\in\hat{B}}\hat{\P}^y(\hat{\xi}(\hat{T})=x,\hat{T}<\infty,\hat{\tau}<\infty)\\
			&=\sum_{x\in\hat{B}}\hat{\P}^y(\hat{\tau}<\infty\,|\,\hat{\xi}(\hat{T})=x,\hat{T}<\infty)\,
			\hat{\P}^y(\hat{\xi}(\hat{T})=x,\hat{T}<\infty)\\
			&=\sum_{x\in\hat{B}}\hat{\P}^x(\hat{\tau}<\infty)\,\hat{\P}^y(\hat{\xi}(\hat{T})=x,\hat{T}<\infty)\\
			&=\sum_{x\in\hat{B}}\hat{\nu}(x)\,\hat{\P}^y(\hat{\xi}(\hat{T})=x,\hat{T}<\infty)
			\geq \epsilon\,\hat{\P}^y(\hat{T}<\infty)
			\geq \epsilon,
		\end{aligned}
	\end{equation}
	which contradicts \eqref{eqn: zero-infimum}. 
\end{proof}

\begin{corollary}{\bf [Equivalence of absorption probabilities]}
	\label{cor: equiv-of-hat-and-orig-sys}
	For any $\eta\in\mathcal{S}$, $\nu(\eta) = 1$ if $\hat{\nu}(\eta)=1$. Furthermore, if \eqref{as: uniformly-bounded-nu-hat-probability} holds, then the converse is true as well.
\end{corollary}

\begin{proof}
	The claim follows from Corollary~\ref{corollary: nu > nu_hat} and Theorem~\ref{theorem: nu_hat > nu}.
\end{proof}


\subsection{Comparison with non-interacting dual}
\label{ss.compindep}

{The goal of this section is to reduce the absorption analysis of $\xi$ and $\hat{\xi}$ in the previous section to equivalent statements involving the \emph{independent RW1} introduced in Definition~\ref{definition: independent RW}. We follow the same comparison method used earlier.}

\begin{theorem}{\bf [Comparison of absorption probabilities]}
	\label{theorem: nu* > nu_hat}
	Let $\nu^*\colon\,\mathcal{S}\to[0,1]$ and $\hat\nu\colon\,\mathcal{S}\to[0,1]$ be defined by
	\begin{equation}
		\nu^*(\eta) := \P^*_\eta(\tau^*<\infty),\quad\hat{\nu}(\eta):=\hat{\P}^\eta(\hat{\tau}<\infty).
	\end{equation}
	Assume that
	\begin{equation}
		\label{as: uniformly-bounded-nu*-probability}
		\begin{aligned}
			\inf\{\nu^*([(i,1),(i,1)])\colon\,i\in\Z^d\}>0.
		\end{aligned}
	\end{equation}
	{For all $\eta\in\mathcal{S}$, if $\hat{\nu}(\eta) = 1$, then $\nu^*(\eta)=1$.}
\end{theorem}

\begin{proof}
	The proof follows a similar argument as in the proof of Theorem~\ref{theorem: nu_hat > nu}. Suppose that $\hat{\nu}(\eta)=1$ and $\nu^*(\eta)<1$. By the strong Markov property,
	\begin{equation}
		\label{eqn: nu*-zero-infimum}
		\begin{aligned}
			\inf_{y\in\S}\nu^*(y)=0.
		\end{aligned}
	\end{equation}
	Since, by Remark~\ref{remark: irreducible hat system}, the process $\hat{\xi}$ started from $\eta$ can visit any configuration $y\in\S$ in finite time with positive probability, we have 
	\begin{equation}
		\label{eqn: constant-nu-hat}
		\begin{aligned}
			\hat{\nu}(y) = 1 \qquad \forall\,y\in\S.
		\end{aligned}
	\end{equation}
	We will show that \eqref{eqn: nu*-zero-infimum} and \eqref{eqn: constant-nu-hat} are contradictory. 
	
	Let $\bar{B}\subset\S$ be defined as
	\begin{equation}
		\begin{aligned}
			\bar{B}:=\Big\{[(i,\alpha),(i,\beta)]\in\S\colon\,\alpha\neq \beta,\nu^*([(i,1),(i,1)])<\nu^*([(i,1),(i,0)])\Big\}\cup\{\circledast\}.
		\end{aligned}
	\end{equation}
	By symmetry and a first-{jump} analysis, we have
	\begin{equation}
		\label{eqn: constant-nu*-on-interchange}
		\begin{aligned}
			\nu^*([(i,1),(i,0)])=\nu^*([(i,0),(i,1)])=\nu^*([(i,0),(i,0)]) \qquad \forall\,i\in\Z^d.
		\end{aligned}
	\end{equation} 
	Let $\hat{T}_{\bar{B}}$ denote the first hitting time of the set $\bar{B}$ for the process $\hat{\xi}$, and let 
	\begin{equation}
		\label{eqn: epsilon-defn-nu*}
		\begin{aligned}
			\bar{\epsilon}:=\inf\{\nu^*(y)\colon\,y\in \bar{B}\}.
		\end{aligned}
	\end{equation}
	By \eqref{as: uniformly-bounded-nu*-probability} and \eqref{eqn: constant-nu*-on-interchange}, $\bar{\epsilon}>0$. Note that if $\hat{Q}$ and $Q^*$ are the generators of the processes $\hat{\xi}$ and $\xi^*$, respectively, then
	\begin{equation}
		\begin{aligned}
			&((\hat{Q}-Q^*)\nu^*)(x)\\ 
			&=\begin{cases}
				\tfrac{\lambda}{M_i}\delta_{i,j}(1-\delta_{\alpha,\beta})[\nu^*([(i,1),(i,0)])-\nu^*([(i,1),(i,1)])], 
				&x=[(i,\alpha),(j,\beta)]\neq\circledast,\\
				0, 
				&\text{ otherwise},
			\end{cases}
		\end{aligned}
	\end{equation}
	where \eqref{eqn: constant-nu*-on-interchange} is used. Moreover, the right-hand side of the above equation is negative whenever $x\notin \bar{B}$. Since $Q^*\nu^*\equiv 0$, we have
	\begin{equation}
		\label{eqn: Q-hat-negative-outside-B}
		\begin{aligned}
			(\hat{Q}\nu^*)(x) \leq 0, \qquad x\notin \bar{B}.
		\end{aligned}
	\end{equation}
	Let $y\in\S$ be fixed arbitrarily, and let the process $\hat{\xi}$ be started from $y$. Since $\nu^*$ is bounded and $\hat{\xi}$ is non-explosive, the process $(M_t)_{t\geq 0}$ with $ M_t := \nu^*(\hat{\xi}(t)) - \int_{0}^{t} (\hat{Q}\nu^*)(\hat{\xi}(s))\,\d s $ is a martingale under the law $\hat{\P}^y$ w.r.t.\ the natural filtration associated to the process $\hat{\xi}$. Hence the stopped process $(M_{t\wedge \hat{T}_{\bar{B}}})_{t\geq 0}$ is also a martingale. By \eqref{eqn: Q-hat-negative-outside-B}, we have $\int_{0}^{t\wedge \hat{T}_{\bar{B}}} (\hat{Q}\nu^*)(\hat{\xi}(s))\,\d s\leq 0$ a.s.\ for any $t\geq 0$. Hence $M_{t\wedge \hat{T}_{\bar{B}}} \geq \nu^*(\hat{\xi}({t\wedge \hat{T}_{\bar{B}}}))$ for any $t\geq 0$. By the martingale property, for any $t> 0$,
	\begin{equation}
		\begin{aligned}
			\nu^*(y) &=\hat{\E}^y[\nu^*(\hat{\xi}(0))] = \hat{\E}^y[M_{t\wedge \hat{T}_{\bar{B}}}]
			\geq \hat{\E}^y[\nu^*(\hat{\xi}({t\wedge \hat{T}_{\bar{B}}}))]\\
			&\geq \hat{\E}^y[\nu^*(\hat{\xi}({\hat{T}_{\bar{B}}}))\mathbf{1}_{\hat{T}_{\bar{B}}< t}]
			\geq \bar{\epsilon}\,\hat{\P}^y(\hat{T}_{\bar{B}}<t)\geq \bar{\epsilon}\, \hat{\P}^y(\hat{\tau}<t),
		\end{aligned}
	\end{equation}
	where in the last inequality we use that $\hat{T}_{\bar{B}}\leq \hat{\tau}$ a.s. Letting $t\to\infty$, we find with the help of \eqref{eqn: constant-nu-hat} that $\nu^*(y)\geq \bar{\epsilon}\,\hat{\P}^y(\hat{\tau}<\infty)=\bar{\epsilon}\,\hat{\nu}(y)=\bar{\epsilon}$, which contradicts \eqref{eqn: nu*-zero-infimum}.
\end{proof}

\begin{theorem}{\bf [Comparison of absorption probabilities]}
	\label{theorem: nu-hat > nu*}
	Let $\nu^*$, $\nu$, $\hat{\nu}$ be the absorption probability of $\xi^*$, $\hat{\xi}$, $\xi$, respectively, i.e.,
	\begin{equation}
		\nu^*(\eta) := \P^*_\eta(\tau^*<\infty),\quad\hat{\nu}(\eta):=\hat{\P}^\eta(\hat{\tau}<\infty),
		\quad\nu(\eta):=\P^\eta(\tau<\infty).
	\end{equation}
	Assume that
	\begin{equation}
		\label{as: strongly-bounded-nu-hat-probability}
		\begin{aligned}
			\inf\{\hat{\nu}([(i,1),(i,0)])\colon\,i\in\Z^d\}>0.
		\end{aligned}
	\end{equation}
	{For all $\eta\in\mathcal{S}$, if $\nu^*(\eta) = 1$, then $\hat{\nu}(\eta)=1$,} and hence $\nu(\eta) = 1$ as well.
\end{theorem}

\begin{proof}
	By Corollary~\ref{corollary: nu > nu_hat}, it suffices to prove that $\hat{\nu}(\eta)=1$. Suppose that this fails. Then, by the strong Markov property,
	\begin{equation}
		\label{eqn: nu-hat-zero-infimum}
		\begin{aligned}
			\inf_{y\in\S}\hat{\nu}(y)=0.
		\end{aligned}
	\end{equation}
	Moreover, since the process $\xi^*$ started from $\eta$ can visit any configuration $y\in\S$ in finite time with positive probability, we have
	\begin{equation}
		\label{eqn: constant-nu*}
		\begin{aligned}
			\nu^*(y) = 1 \qquad \forall\, y\in\S.
		\end{aligned}
	\end{equation}
	We will show that \eqref{eqn: nu-hat-zero-infimum} and \eqref{eqn: constant-nu*} are contradictory. 
	
	Let $B^\prime\subset\S$ be defined as
	\begin{equation}
		\begin{aligned}
			B^\prime:=\Big\{[(i,\alpha),(i,\beta)]\in\S\colon\,\alpha\neq \beta,\hat{\nu}([(i,1),(i,1)])
			\geq\hat{\nu}([(i,1),(i,0)])\Big\}\cup\{\circledast\}.
		\end{aligned}
	\end{equation}
	By symmetry and a first-{jump} analysis, we have
	\begin{equation}
		\label{eqn: constant-nu-hat-on-interchange}
		\begin{aligned}
			\hat{\nu}([(i,1),(i,0)])=\hat{\nu}([(i,0),(i,1)])=\hat{\nu}([(i,0),(i,0)]) \qquad \forall\,i\in\Z^d.
		\end{aligned}
	\end{equation} 
	Let $T^*_{B^\prime}$ denote the first hitting time of the set $B^\prime$ for the process $\xi^*$, and let 
	\begin{equation}
		\label{eqn: epsilon-defn-nu-hat}
		\begin{aligned}
			\epsilon^\prime:=\inf\{\hat{\nu}(y)\colon\,y\in B^\prime\}.
		\end{aligned}
	\end{equation}
	By \eqref{as: strongly-bounded-nu-hat-probability} and \eqref{eqn: constant-nu-hat-on-interchange}, we have $\epsilon^\prime>0$. Note that if $\hat{Q}$ and $Q^*$ are the generators of the processes $\hat{\xi}$ and $\xi^*$, respectively, then
	\begin{equation}
		\begin{aligned}
			&((Q^*-\hat{Q})\hat{\nu})(x)\\ 
			&=
			\begin{cases}
				\tfrac{\lambda}{M_i}\delta_{i,j}(1-\delta_{\alpha,\beta})[\hat{\nu}([(i,1),(i,1)])-\hat{\nu}([(i,1),(i,0)])],
				&x=[(i,\alpha),(j,\beta)]\neq\circledast,\\
				0,&\text{ otherwise},
			\end{cases}
		\end{aligned}
	\end{equation}
	where we use \eqref{eqn: constant-nu-hat-on-interchange}. Moreover, the right-hand side side of the above equation is negative whenever $x\notin B^\prime$. Since $\hat{Q}\hat{\nu}\equiv 0$, we have
	\begin{equation}
		\label{eqn: Q*-negative-outside-B}
		\begin{aligned}
			(Q^*\hat{\nu})(x) \leq 0, \qquad x\notin B^\prime.
		\end{aligned}
	\end{equation}
	Let $y\in\S$ be fixed arbitrarily, and let the process $\xi^*$ be started from $y$. Since $\hat{\nu}$ is bounded and $\xi^*$ is non-explosive, the process $(M_t)_{t\geq 0}$ with \begin{equation}
		M_t := \hat{\nu}(\xi^*(t)) - \int_{0}^{t} (Q^*\hat{\nu})(\xi^*(s))\,\d s
	\end{equation}
	is a martingale under the law $\P^*_y$ w.r.t.\ the natural filtration associated to the process $\xi^*$. Hence the stopped process $(M_{t\wedge T^*_{B^\prime}})_{t\geq 0}$ is also a martingale. By \eqref{eqn: Q*-negative-outside-B}, we have $\int_{0}^{t\wedge T^*_{B^\prime}} (Q^*\hat{\nu})(\xi^*(s))\,\d s\leq 0$ a.s.\ for any $t\geq 0$. Hence $M_{t\wedge T^*_{B^\prime}} \geq \hat{\nu}(\xi^*({t\wedge T^*_{B^\prime}}))$ for any $t\geq 0$. By the martingale property, for any $t> 0$,
	\begin{equation}
		\begin{aligned}
			\hat{\nu}(y) &=\E^*_y[\hat{\nu}(\xi^*(0))] = \E^*_y[M_{t\wedge T^*_{B^\prime}}]
			\geq \E^*_y[\hat{\nu}(\xi^*({t\wedge T^*_{B^\prime}}))]\\
			&\geq \E^*_y[\hat{\nu}(\xi^*({T^*_{B^\prime}}))\mathbf{1}_{T^*_{B^\prime}< t}]
			\geq \epsilon^\prime\,\P^*_y(T^*_{B^\prime}<t)\geq \epsilon^\prime\, \P^*_y(\tau^*<t),
		\end{aligned}
	\end{equation}
	where in the last inequality we use that $T^*_{B^\prime}\leq \tau^*$ a.s. Letting $t\to\infty$, we find via \eqref{eqn: constant-nu*} that $\hat{\nu}(y)\geq \epsilon^\prime\,\P^*_y(\tau^*<\infty)=\epsilon^\prime\,\nu^*(y)=\epsilon^\prime$, which contradicts \eqref{eqn: nu-hat-zero-infimum}. 
\end{proof}

\begin{remark}
	{\rm Theorem~\ref{theorem: nu-hat > nu*} tells us that coalescence of independent particles is sufficient for coalescence of interacting particles. The condition {in \eqref{as: strongly-bounded-nu-hat-probability}} is stronger, because it requires control on the growth of both $N_i$ and $M_i$.} \hfill$\Box$
\end{remark}


\subsection{Conclusion}
\label{ss.concl}

\begin{theorem}{\bf [Equivalence of absorption probabilities]} 
	\label{theorem: necessity of coalescence of independent particles}
	Let $\nu^*$, $\nu$ and $\hat{\nu}$ be the functions defined by
	\begin{equation}
		\nu^*(\eta) := \P^*_\eta(\tau^*<\infty),\quad\hat{\nu}(\eta):=\hat{\P}^\eta(\hat{\tau}<\infty),\quad\nu(\eta):=\P^\eta(\tau<\infty).
	\end{equation}
	If
	\begin{itemize}
		\item[{\rm (a)}]
		$\inf\{\hat{\nu}([(i,1),(i,1)])\colon\,i\in\Z^d\}>0$,
		\item[{\rm (b)}]  
		$\inf\{\nu^*([(i,1),(i,1)])\colon\,i\in\Z^d\}>0$,
	\end{itemize}
	then $\nu^*(\eta)=1$ whenever $\nu(\eta) = 1$ for some $\eta\in\mathcal{S}$. If $\inf\{\hat{\nu}([(i,1),(i,0)])\colon\,i\in\Z^d\}>0$, then the converse is true as well. 
\end{theorem}

\begin{proof}
	The forward direction follows by combining Theorem~\ref{theorem: nu_hat > nu} and Theorem~\ref{theorem: nu* > nu_hat}. The reverse direction is a direct consequence of Theorem~\ref{theorem: nu-hat > nu*} and Corollary~\ref{corollary: nu > nu_hat}.
\end{proof}

\begin{remark}
	\label{rem:weaker}
	{\rm Theorem~\ref{theorem: necessity of coalescence of independent particles} tells us that if the interacting particle system coalesces with probability $1$, then it is necessary that two independent particles coalesce with probability $1$. The {first} two conditions are trivially satisfied when $\sup_{i\in\Z^d}  N_i<\infty$. If, furthermore, $\sup_{i\in\Z^d} M_i<\infty$, then the third condition is satisfied as well.} \hfill$\Box$
\end{remark}

We conclude this section by providing conditions on the sizes of the active and the dormant populations that are weaker than the ones mentioned in Remark~\ref{rem:weaker}, and under which the assumptions in Theorem~\ref{theorem: necessity of coalescence of independent particles} are satisfied.

\begin{theorem}{\bf [Lower bound on absorption probabilities]} 
	\label{thm:lower-bound-on-probability}
	Let $\hat{\nu}$ and $\nu^*$ be the functions defined by
	\begin{equation}
		\hat{\nu}(\eta):=\hat{\P}^\eta(\hat{\tau}<\infty),\quad\nu^*(\eta) := \P^*_\eta(\tau^*<\infty).
	\end{equation}
	If the sizes of the active populations $(N_i)_{i\in\Z^d}$ are non-clumping, i.e.,
	\begin{equation}
		\inf_{i\in\Z^d} \sum_{\|j-i\|\leq R} \tfrac{1}{N_j}>0 \text{ for some } R<\infty,
	\end{equation}
	then
	\begin{itemize}
		\item[\rm{(a)}] 
		$\inf\{\hat{\nu}([(i,1),(i,1)])\colon\,i\in\Z^d\}>0$.
		\item[\rm{(b)}] 
		$\inf\{\nu^*([(i,1),(i,1)])\colon\,i\in\Z^d\}>0$.
	\end{itemize}
	Furthermore, if the relative strengths of the seed-banks are bounded, i.e.,  
	\begin{equation}
		\sup_{i\in\Z^d}\frac{M_i}{N_i}<\infty,
	\end{equation}
	then
	\begin{itemize}
		\item[\rm{(i)}] 
		$\inf\{\hat{\nu}([(i,1),(i,0)])\colon\,i\in\Z^d\}>0$.
		\item[\rm{(ii)}] $\inf\{\nu^*([(i,1),(i,0)])\colon\,i\in\Z^d\}>0$.
	\end{itemize}
\end{theorem}

Before we give the proof of Theorem~\ref{thm:lower-bound-on-probability} we derive a series representation of the absorption probabilities $\nu^*$ and $\hat{\nu}$ of the respective processes $\xi^*$ and $\hat{\xi}$.

\begin{lemma}{\bf [Series representation]}
	\label{lemma:absorption-probability-repr}
	Let $\nu^*$ and $\hat{\nu}$ be the functions defined by
	\begin{equation}
		\nu^*(\eta) := \P^*_\eta(\tau^*<\infty),\quad\hat{\nu}(\eta):=\hat{\P}^\eta(\hat{\tau}<\infty).
	\end{equation}
	For $i\in\Z^d$, let $R_i^*$ (respectively, $\hat{R}_i$) be the total number of visits to the state $[(i,1),(i,1)]\in\mathcal{S}$ made by the jump chain associated to the process $\xi^*$ (respectively, $\hat{\xi}$). Then, for $\eta\in\mathcal{S}\backslash\{\circledast\}$,
	\begin{itemize}
		\item[\rm{(a)}] 
		$\displaystyle\nu^*(\eta) = \sum_{i\in\Z^d}\tfrac{1}{2(c+\lambda)N_i+1}\,\E_\eta^*[R_i^*]$.
		\item[\rm{(b)}]$\displaystyle\hat{\nu}(\eta) = \sum_{i\in\Z^d}\tfrac{1}{2(c+\lambda)N_i+1}\,\hat{\E}^\eta[\hat{R}_i]$,
	\end{itemize}
	where $c$ is the total migration rate defined in Assumption~\ref{assumpt1}, and expectations are taken w.r.t.\ the respective laws of the jump chains associated to the processes $\xi^*$ and $\hat{\xi}$.
\end{lemma}

\begin{proof}
	We only prove part (a), because the proof of part (b) is the same. Let $\eta\in\mathcal{S}\backslash\{\circledast\}$ be fixed, and let $X^* := (X_n^*)_{n\in\N_0}$ be the embedded jump chain associated to the process $\xi^*$ started at state $\eta$. Since $X^*$ is absorbed to $\circledast$ if and only if $\xi^*$ is absorbed, it suffices to analyse $X^*$. Let $T:=\inf\{n\in\N_0\colon \,X_n^*=\circledast\}$ be the absorption time of $X^*$. Note that, because the absorbing state $\circledast$ can be reached in one step {only} from the states $\{[(i,1),(i,1)]\colon\,i\in\Z^d\}\subset\mathcal{S}$, for all $n\in\N$ we have
	\begin{equation}
		\begin{aligned}
			\P_\eta^*(T=n) &= \sum_{i\in\Z^d}\P_\eta^*(X^*_{n-1}=[(i,1),(i,1)], T=n)\\
			&=\sum_{i\in\Z^d}\P_\eta^*(X^*_n=\circledast\,|\,X^*_{n-1}=[(i,1),(i,1)])\,\P_\eta^*( X_{n-1}^*=[(i,1),(i,1)])\\
			&=\sum_{i\in\Z^d}\tfrac{1}{2(c+\lambda)N_i+1}\,\P_\eta^*( X_{n-1}^*=[(i,1),(i,1)]),
		\end{aligned}
	\end{equation}
	where in the last equality we use that, by the Markov property, 
	\begin{equation}
		\P_\eta^*(X^*_n=\circledast\,|\,X^*_{n-1}=[(i,1),(i,1)])=\P_{[(i,1),(i,1)]}^*(X^*_1=\circledast)=\frac{1}{2(c+\lambda)N_i+1}.
	\end{equation}
	Using that $\eta\neq\circledast$, we get
	\begin{equation}
		\begin{aligned}
			\nu^*(\eta)
			&=\P^*_\eta(T<\infty)=\sum_{n\in\N} \P_\eta^*(T=n)\\
			&=\sum_{n\in\N} \sum_{i\in\Z^d}\tfrac{1}{2(c+\lambda)N_i+1}\,\P_\eta^*( X_{n-1}^*=[(i,1),(i,1)])\\
			&=\sum_{i\in\Z^d}\tfrac{1}{2(c+\lambda)N_i+1}\sum_{n\in\N} \P_\eta^*( X_{n-1}^*=[(i,1),(i,1)])
			=\sum_{i\in\Z^d}\tfrac{1}{2(c+\lambda)N_i+1}\,\E^*_\eta[R_i^*],
		\end{aligned}
	\end{equation}
	where in the fourth equality we interchange the two sums using Fubini's theorem, and in the last equality we use
	\begin{equation}
		\E_\eta^*[R_i^*] = \sum_{n\in\N_0}\P_\eta^*(X^*_n=[(i,1),(i,1)]),\quad\quad i\in\Z^d.
	\end{equation}
\end{proof}

\begin{proof}[Proof of Theorem~\ref{thm:lower-bound-on-probability}]
	We only prove parts (a) and (i), because the proof of parts (b) and (ii) is the same. Let $\hat{X}:=(\hat{X}_n)_{n\in\N_0}$ be the embedded jump chain associated to the process $\hat{\xi}$. For $j\in\Z^d$, let $\hat{R}_j$ be the total number of visits made by $\hat{X}$ to the state $[(j,1),(j,1)]$. We first show that, for any $i,j\in\Z^d$,
	\begin{equation}
		\hat{\E}^{[(i,1),(i,1)]}[\hat{R}_j]\geq \sum_{n\in\N} m^{2n}a_n(i,j)^2,
	\end{equation}
	where $m:=\tfrac{c}{2(c+\lambda)+1}$. Note that, in the process $\hat{\xi}$, each of the two particles moves from $i$ to $j$ at rate $a(i,j)$ while in the active state, and becomes dormant at rate $\lambda$ when the two particles are not on top of each other with one active and the other dormant. Thus, for $i,j,k\in\Z^d$ and $n\in\N$,
	\begin{equation}
		\begin{aligned}
			&\hat{\P}^{[(k,1),(i,1)]}(\hat{X}_n=[(k,1),(j,1)])\\
			&\geq\sum_{l\neq i}\hat{\P}^{[(k,1),(i,1)]}(\hat{X}_1=[(k,1),(l,1)])\,\hat{\P}^{[(k,1),(l,1)]}(\hat{X}_{n-1}=[(k,1),(j,1)])\\
			&=\sum_{l\neq i}\tfrac{c}{2(c+\lambda)+(1/N_i)\delta_{k,i}}\tfrac{a(i,l)}{c}\,\hat{\P}^{[(k,1),(l,1)]}(\hat{X}_{n-1}=[(k,1),(j,1)])\\
			&\geq m\sum_{l\neq i}a_1(i,l)\,\hat{\P}^{[(k,1),(l,1)]}(\hat{X}_{n-1}=[(k,1),(j,1)]),
		\end{aligned}
	\end{equation}
	where $a_1(\cdot\,,\cdot):=\tfrac{a(\cdot\,,\,\cdot)}{c}$ is the transition kernel of the embedded chain associated to the continuous-time random walk on $\Z^d$ with rates $a(\cdot\,,\cdot)$. Using the above recursively, we obtain that, for any $i,j,k\in\Z^d$ and $n\in\N$,
	\begin{equation}
		\hat{\P}^{[(k,1),(i,1)]}(\hat{X}_n=[(k,1),(j,1)])\geq m^n\,a_n(i,j).
	\end{equation}
	Therefore, applying the above twice, for $i,j\in\Z^d$ we have
	\begin{equation}
		\begin{aligned}
			\hat{\P}^{[(i,1),(i,1)]}(\hat{X}_{2n}=[(j,1),(j,1)])
			&\geq \hat{\P}^{[(i,1),(i,1)]}(\hat{X}_{n}=[(i,1),(j,1)])\,\hat{\P}^{[(i,1),(j,1)]}(\hat{X}_{n}=[(j,1),(j,1)])\\
			&\geq m^n a_n(i,j)\,\hat{\P}^{[(j,1),(i,1)]}(\hat{X}_{n}=[(j,1),(j,1)])
			\geq m^{2n}a_n(i,j)^2.
		\end{aligned}
	\end{equation}
	Hence, for $i,j\in\Z^d$,
	\begin{equation}
		\begin{aligned}
			\hat{\E}^{[(i,1),(i,1)]}[\hat{R}_j] 
			&=  \sum_{n\in\N_0}\hat{\P}^{[(i,1),(i,1)]}(\hat{X}_n=[(j,1),(j,1)])\\
			&\geq\sum_{n\in\N_0}\hat{\P}^{[(i,1),(i,1)]}(\hat{X}_{2n}=[(j,1),(j,1)])
			\geq \sum_{n\in\N} m^{2n}a_n(i,j)^2.
		\end{aligned}
	\end{equation}
	Finally, substituting the above into the series representation of $\hat{\nu}$ in part (b) of Lemma~\ref{lemma:absorption-probability-repr}, we obtain that, for $i\in\Z^d$,
	\begin{equation}
		\label{eqn:lower-bound-on-nu-hat}
		\begin{aligned}
			\hat{\nu}([(i,1),(i,1)])
			&=\sum_{j\in\Z^d}\tfrac{1}{2(c+\lambda)N_j+1}\hat{\E}^{[(i,1),(i,1)]}[\hat{R}_j]\\
			&\geq \sum_{j\in\Z^d}\tfrac{1}{2(c+\lambda)N_j+1} \sum_{n\in\N} m^{2n}a_n(i,j)^2\\
			&\geq \tfrac{1}{2(c+\lambda)+1}\sum_{j\in B_R(i)}\tfrac{1}{N_j} \sum_{n\in\N} m^{2n}a_n(0,j-i)^2
			\geq \epsilon_R\sum_{j\in B_R(i)}\tfrac{1}{N_j},
		\end{aligned}
	\end{equation}
	where {$B_R(i) := \{j\in\Z^d\colon\,\|j-i\|\leq R\}$ and}
	\begin{equation}
		\epsilon_R:=\min\Big\{\tfrac{1}{2(c+\lambda)+1}\sum_{n\in\N} m^{2n}a_n(0,l)^2\colon\,l\in B_R(0)\Big\}>0.
	\end{equation}
	Since, by assumption, $(N_i)_{i\in\Z^d}$ are non-clumping, the right-hand side of \eqref{eqn:lower-bound-on-nu-hat} is bounded away from zero irrespective of the choice $i\in\Z^d$, and so part (a) is proved. 
	
	To prove part (i), by doing a first-{jump} analysis of the process $\hat{X}$ we get that, for $i\in\Z^d$, 
	\begin{equation}
		\hat{\nu}([(i,1),(i,0)])\geq \hat{\P}^{[(i,1),(i,0)]}(\hat{X}_1=[(i,1),(i,1)])\,\hat{\nu}([(i,1),(i,1)])
		=\tfrac{\lambda K_i}{c+\lambda+\lambda K_i}\hat{\nu}([(i,1),(i,1)]),
	\end{equation}
	where $K_i=\tfrac{N_i}{M_i}$. Thus, if $(N_i)_{i\in\Z^d}$ are non-clumping and $\sup_{i\in\Z^d}K_i^{-1}<\infty$, then
	\begin{equation}
		\hat{\nu}([(i,1),(i,0)])\geq \frac{\lambda}{\lambda+(c+\lambda)(\sup_{i\in\Z^d}K_i^{-1})}\inf\{\hat{\nu}([(j,1),(j,1)])\colon\,j\in\Z^d\},
	\end{equation}
	which is bounded away from zero uniformly in $i\in\Z^d$, and so part (i) follows.
\end{proof}


\section{Proofs: clustering criterion and clustering regime}
\label{s.coalbd}

In this section we prove our two main theorems, namely, Theorem~\ref{T.Alt-Dichotomy} and Theorem~\ref{T.Coal_recc} with the help of the results that were obtained in Section~\ref{s.duals} by comparing various auxiliary duals.

\begin{proof}[Proof of Theorem~\ref{T.Alt-Dichotomy}]
	{Note (see Remark~\ref{remark:equivalence-of-absorption})} that the system clusters if and only if the two-particle process $\xi$ defined in Definition~\ref{definition: interacting RW1} is absorbed to $\circledast$ with probability 1. Let $\hat\xi$ be the auxiliary two-particle process defined in Definition~\ref{definition: interacting RW2}, and $\hat{\nu}(\eta)$ (respectively, $\nu(\eta)$) be the absorption probability of the process $\hat{\xi}$ (respectively, $\xi$) started from state $\eta\in G\times G$. The system $Z$ clusters if and only if $\nu(\eta) = 1$ for any state $\eta\in G\times G$. By the forward direction of Corollary~\ref{cor: equiv-of-hat-and-orig-sys}, we have that $\nu(\eta)=1$ whenever $\hat{\nu}(\eta)=1$, and hence the forward direction of Theorem~\ref{T.Alt-Dichotomy} follows. To prove the converse we note that, under the non-clumping assumption of the active populations sizes $(N_i)_{i\in\Z^d}$ in \eqref{eqn:non-clumping criterion}, \eqref{as: uniformly-bounded-nu-hat-probability} in Corollary~\ref{cor: equiv-of-hat-and-orig-sys} holds by part (a) of Theorem~\ref{thm:lower-bound-on-probability}, and hence $\hat{\nu}(\eta)=1$ whenever $\nu(\eta)=1$, so that the converse follows as well.
\end{proof}

\subsection{Independent particle system and clustering regime.}\label{ss.indep-sys}

In order to prove Theorem~\ref{T.Coal_recc}, we need to take a closer look at the non-interacting two-particle process $\xi^*$ introduced in Definition~\ref{definition: independent RW}. In what follows we briefly describe the process $\xi^*$ and derive conditions under which the process $\xi^*$ is absorbed with probability 1.

\begin{figure}[htbp]
	\begin{center}
		\begin{tikzpicture}[xscale=0.8,yscale=0.8]
			\draw (0,-1.5) rectangle (16.5,1.5); 
			\pgflowlevelsynccm
			\draw[very thick,blue] (1.5,0) -- (4.7,0);
			\draw[very thick,red] (4.7,0) -- (5.9,0);
			\draw [decorate,decoration={brace,amplitude=3pt},yshift=-0.15cm]
			(5.9,0) -- node [below,yshift=-0.1cm]
			{$D_1$} (4.7,0);
			\draw[very thick,blue] (5.9,0) -- (8.3,0);
			\draw[very thick,red] (8.3,0) -- (11,0);
			\draw [decorate,decoration={brace,amplitude=3pt},yshift=-0.15cm]
			(11,0) -- node [below,yshift=-0.1cm]
			{$D_2$} (8.3,0);
			\draw[very thick,blue] (11,0) -- (13.7,0);
			\draw[very thick,red] (13.7,0) -- (15.5,0);
			\draw[fill,blue] (1.5,0) circle (1pt) node[above] {$(i_0,1)$};
			\draw[fill,black] (3.7,0) circle (2pt) node[above] {$(i_1,1)$};
			\draw[fill,black] (6.2,0) circle (2pt) node[above] {$(i_2,1)$};
			\draw[fill,black] (8,0) circle (2pt) node[above,xshift=-0.2cm] {$(i_3,1)$};
			\draw[fill,black] (11.5,0) circle (2pt) node[above] {$(i_4,1)$};
			\draw[fill,black] (13,0) circle (2pt) node[above] {$(i_5,1)$};
			\draw [decorate,decoration={brace,amplitude=3pt},yshift=-0.15cm]
			(15.5,0) -- node [below,yshift=-0.1cm] {$D_3$} (13.7,0);
			\draw[dotted] (1.5,0) -- (1.5,-1)  node[below] {time $0$};
			\draw[dotted] (14.6,0.5) -- (14.6,-1)  node[below] {time $t$};
		\end{tikzpicture}
	\end{center}
	\caption{\small Evolution of a single particle started at location $i_0$ in the active state. Red and blue lines denote the dormant and the active phases of the particle. Each dot represents a migration step.}
	\label{fig:1}
\end{figure}
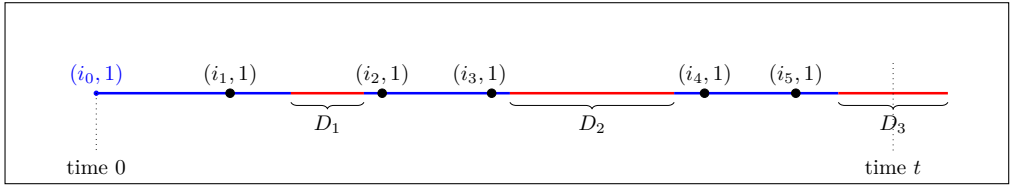

We recall from Definition~\ref{definition: independent RW} that the process $\xi^* = (\xi^*(t))_{t\geq 0}$ is a continuous-time Markov process on the state space $\mathcal{S} = (G\times G)\cup\{\circledast\}$ with $G=\Z^d\times\{0,1\}$. Here, $\xi^*(t)=[(i,\alpha),(j,\beta)]$ captures the location ($i,j \in \Z^d$) and the state ($\alpha,\beta\in\{0,1\}$) of two independent particles at time $t$, where $0$ stands for dormant state and $1$ stands for active state, respectively. The evolution of the two independent particles is governed by the following transitions (see Fig.~\ref{fig:1}):

\begin{itemize}
	\item\textbf{(Migration)} 
	Each particle migrates from location $i$ to $j$ at rate $a(i,j)$ while being active.
	\item \textbf{(Active to Dormant)} 
	An active particle becomes dormant (without changing location) at rate $\lambda$.
	\item \textbf{(Dormant to Active)} 
	A dormant particle at location $i$ becomes active (without changing location) at rate $\lambda K_i$.
	\item \textbf{(Coalescence)} 
	The two particles coalesce with each other, and are absorbed to the state $\circledast$, at rate $\tfrac{1}{N_i}$ when they are both at location $i$ and both active.
\end{itemize}

The following lemma tells that if the mean wake-up time of a dormant particle is uniformly bounded over all the locations in $\Z^d$, then the accumulated activity time of a single particle increases linearly in time.

\begin{lemma}{\bf[Linear activity time]}
	\label{lemma:linear-activity-time}
	Let $S(t)$ be the total accumulated time spent in the active state during the time interval $[0,t]$ by a single particle that evolves according to the first three transitions described above. If $\sup_{i\in\Z^d}K_i^{-1}<\infty$, then
	\begin{equation}
		\liminf\limits_{t\to\infty}\frac{S(t)}{t}{\geq\frac{1}{1+K^{-1}}}\quad\text{a.s.},
	\end{equation}
	{where $K^{-1}:=\sup_{i\in\Z^d}K_i^{-1}$.}
\end{lemma}

\begin{proof}
	We prove the claim with the help of coupling in combination with a renewal argument. Let $(T_n)_{n\in\N}$ and $(D_n)_{n\in\N}$ be the successive time periods during which the particle is in the active and the dormant state, respectively (see Fig.~\ref{fig:1}). Note that $(T_n)_{n\in\N}$ are i.i.d.\ exponential random variables with mean $\frac{1}{\lambda}$. Also note that $D_n$ is exponentially distributed with $\E[D_n]\leq (\lambda K)^{-1}$, because the particle wakes up from the dormant state at rate {$\lambda K_i\geq \lambda K$} when it is at location $i$. Hence, using monotone coupling of exponential random variables, we can construct a sequence $(U_n)_{n\in\N}$ of i.i.d.\ exponential random variables on the same probability space with mean {$(\lambda K)^{-1}$} such that $D_n\leq U_n$ a.s.\ for all $n\in\N$. Consider the alternating renewal process $(R_t)_{t\geq 0}$ that takes value 0 (respectively, 1) during the time intervals $(T_n)_{n\in\N}$ (respectively, $(U_n)_{n\in\N}$), and let $D(t):=t-S(t)$ be the total accumulated time spent in the dormant state during the time interval $[0,t]$. Note that, because $D_n\leq U_n$ a.s.\ for $n\in\N$, we have
	\begin{equation}
		D(t)\leq\int_{0}^{t}\mathbf{1}_{\{R_s=1\}}\,\d s.
	\end{equation} 
	By applying the renewal reward theorem (see e.g. \cite[Section 2b, Chapter VI]{Asm03} or \cite[Theorem 1, Section 10.5]{GS001}) to the process $(R_t)_{t\geq 0}$, we see that
	\begin{equation}
		\limsup_{t\to\infty}\frac{D(t)}{t}\leq \lim\limits_{t\to\infty}\frac{1}{t}\int_{0}^{t}\mathbf{1}_{\{R_s=1\}}\,\d s 
		= \frac{\E[U_n]}{\E[T_n]+\E[U_n]} = {\frac{\tfrac{1}{\lambda K}}{\tfrac{1}{\lambda}+\tfrac{1}{\lambda K}} 
			= \frac{1}{1+K}} \qquad\text{a.s.}
	\end{equation}
	Hence
	\begin{equation}
		\liminf_{t\to\infty}\frac{S(t)}{t}= 1-\limsup_{t\to\infty}\frac{D(t)}{t}\geq{\frac{1}{1+K^{-1}}}>0 \qquad \text{a.s.}
	\end{equation}
\end{proof}

{Before we proceed with the proof of Theorem~\ref{T.Coal_recc}, we need the following lemma, which roughly tells that under the same assumption as in Lemma~\ref{lemma:linear-activity-time} and under Assumption~\ref{assumpt2}, the presence of dormancy does not affect the recurrence behaviour of a single particle evolving according to the symmetrised migration kernel.}

{
	\begin{lemma}{\bf[Recurrence]}
		\label{lemma:recurrence-of-single-particle}
		Let $S(t)$ be the total accumulated time spent in the active state during the time interval $[0,t]$ by a single particle that evolves according to the first three transitions of the independent particle system described earlier, with migration controlled by the symmetrised kernel $\hat{a}(\cdot\,,\cdot)$. If $K^{-1}<\infty$ and Assumption~\ref{assumpt2} holds, then
		\begin{equation}
			\E\left[\int_0^\infty\hat{a}_{S(t)}(0,0)\,\d t\right] = \infty\text{ if and only if } \int_0^\infty\hat{a}_{t}(0,0)\,\d t=\infty,
		\end{equation}
		where the expectation is taken w.r.t.\ the law of the process describing the evolution of the particle.
	\end{lemma}
}

{\begin{proof}
		We prove the stronger statement that, for some constants $C_1,C_2>0$,
		\begin{equation}
			C_1\leq \liminf_{t\to\infty}\frac{\hat{a}_{S(t)}(0,0)}{\hat{a}_t(0,0)} \leq 
			\limsup_{t\to\infty}\frac{\hat{a}_{S(t)}(0,0)}{\hat{a}_t(0,0)}\leq C_2\quad \text{a.s.,}
		\end{equation}
		from which the claim follows. Let $\delta:=\frac{1}{1+K^{-1}}  \in (0,1)$. By Assumption~\ref{assumpt2}, we have
		\begin{equation}
			\label{eqn:regular-variation-uniform-convergence}
			\lim\limits_{t\to\infty}\frac{\hat{a}_{pt}(0,0)}{\hat{a}_{t}(0,0)} =\frac{1}{p^{\sigma}},
		\end{equation}
		where the convergence is uniform in $p\in[\tfrac{\delta}{2},1]$ (see e.g., \cite[Theorem 1.5.2, Section 1.5]{Bingham1987}). Thus, we can find a $T>0$ such that, for all $t\geq T,$
		\begin{equation}
			\sup_{p\in[\tfrac{\delta}{2},1]}\left|\frac{\hat{a}_{pt}(0,0)}{\hat{a}_{t}(0,0)} -p^{-\sigma}\right|<\frac{1}{2}.
		\end{equation}
		In particular, for all $t\geq T$ and $p\in[\tfrac{\delta}{2},1]$,
		\begin{equation}
			\label{eqn:sandwich-inequality}
			\frac{1}{2} \leq \frac{\hat{a}_{pt}(0,0)}{\hat{a}_t(0,0)}\leq \Big(\frac{2}{\delta}\Big)^\sigma+\frac{1}{2}.
		\end{equation}
		Since, by Lemma~\ref{lemma:linear-activity-time}, $\liminf_{t\to\infty}\tfrac{S(t)}{t}\geq \delta$ a.s., we have that $\tfrac{S(t)}{t}\in[\tfrac{\delta}{2},1]$ eventually a.s.\ as $t\to\infty$. Combining this with \eqref{eqn:sandwich-inequality}, we obtain
		\begin{equation}
			\liminf_{t\to\infty} \frac{\hat{a}_{S(t)}(0,0)}{\hat{a}_t(0,0)} = \liminf_{t\to\infty} \frac{\hat{a}_{(S(t)/t)\, t}(0,0)}{\hat{a}_t(0,0)} \geq \frac{1}{2},\quad\text{a.s.,}
		\end{equation}
		and similarly $\limsup_{t\to\infty}\frac{\hat{a}_{S(t)}(0,0)}{\hat{a}_t(0,0)}\leq \big(\tfrac{2}{\delta}\big)^\sigma+\frac{1}{2}$ a.s.
	\end{proof}
}

{
	\begin{remark}
		\label{remark:recurrence-of-sym-kernel}
		{\rm The proof of the above lemma only uses the regular variation of $\hat{a}_t(0,0)$ at infinity and the fact that $\liminf_{t\to\infty}\tfrac{S(t)}{t}>\delta$ a.s.\ for some $\delta\in(0,1)$. Thus, if $S^\prime(\cdot)$ is an independent copy of $S(\cdot)$, then we also have that
			\begin{equation}
				\E\left[\int_0^\infty\hat{a}_{S(t)+S^\prime(t)}(0,0)\,\d t\right] = \infty\text{ if and only if } \int_0^\infty\hat{a}_{2t}(0,0)\,\d t=\infty,
			\end{equation}
			which is again equivalent to $\hat{a}(\cdot\,,\cdot)$ being recurrent.}\hfill$\Box$
	\end{remark}
}

The following result provides a necessary and sufficient condition for the absorption of the process $\xi^*$.

\begin{theorem}{\bf [Clustering regime]}
	\label{theorem: crit-for-abs-of-indep}
	Suppose that $K^{-1} = \sup_{i\in\Z^d} K_i^{-1}<\infty$ and Assumption~\ref{assumpt2} holds. If the process $\xi^*$ is absorbed to $\circledast$ with probability 1, then it is necessary that the symmetrised kernel $\hat{a}(\cdot\,,\cdot)$ is recurrent, i.e.,
	\begin{equation}
		\label{eqn: recurrent_kernel}
		\int_{0}^{\infty}\hat{a}_t(0,0)\,\d t = \infty.
	\end{equation}
	Furthermore, if $(N_i)_{i\in\Z^d}$ satisfies the non-clumping condition in \eqref{eqn:non-clumping criterion} and $a(\cdot\,,\cdot)$ is symmetric, then \eqref{eqn: recurrent_kernel} is also sufficient.
\end{theorem}

\begin{proof}
	Without loss of generality we may assume that the process starts at the state $\eta := [(0,1),(0,1)]$, i.e., both particles are initially at the origin $0\in\Z^d$ and in the active state. Since the process $\xi^*$ has a positive rate of absorption only when the two independent particles are on top of each other and active, for the absorption probability to be equal to 1 it is necessary that, in the process where coalescence is switched off, the two independent particles meet infinitely often on the same location with probability 1. Let $S(t)$ and $S^\prime(t)$ denote the total accumulated time spent in the active state by the two independent particles (where coalescence is switched off) during the time interval $[0,t]$. Since the two particles move according to $a(\cdot\,,\cdot)$ only when they are active, the total average time during which the two particles are on top of each other is given by
	\begin{equation}
		\label{eqn:average-intersection-time}
		I := \int_{0}^{\infty} f(t)\,\d t,
	\end{equation}
	where $f(t)$ is the probability that the two particles are on the same location at time $t$, which is given by
	\begin{equation}
		\label{eqn:probability-of-meeting}
		f(t):=\E^*_\eta\Big[\sum_{i\in\Z^d}a_{S(t)}(0,i)a_{S^\prime(t)}(0,i)\Big].
	\end{equation}
	Thus, for the process $\xi^*$ to be absorbed with probability 1, it is necessary that $I=\infty$. 
	
	Let us define
	\begin{equation}
		M(t):=S(t)\wedge S^\prime(t),\quad L(t) := [S(t)\vee S^\prime(t)]-[S(t)\wedge S^\prime(t)] = |S(t)-S^\prime(t)|.
	\end{equation}
	Note that
	\begin{equation}
		\label{eqn:probability-of-meeting-1}
		\sum_{i\in\Z^d}a_{S(t)}(0,i)a_{S^\prime(t)}(0,i) = \sum_{i\in\Z^d}\hat{a}_{2M(t)}(0,i)a_{L(t)}(i,0),
	\end{equation}
	because the difference of two continuous-time random walks started at the origin that move independently in $\Z^d$ with rates $a(\cdot\,,\cdot)$ has distribution $\hat{a}_{2M(t)}(0,\cdot)$ at time $M(t)$ (because $a(\cdot\,,\cdot)$ is translation-invariant), and in order for the particle with the largest activity time to meet the other particle at the activity time $S(t)\vee S^\prime(t) = M(t)+L(t)$, it must bridge the difference in the remaining time $L(t)$. We use the Fourier representation of the transition probability kernel $b(\cdot\,,\cdot)$, defined by
	\begin{equation}
		b(i,j) := \frac{a(i,j)}{c}\mathbf{1}_{i\neq j},\qquad i,j\in\Z^d,
	\end{equation} 
	to further simplify the expression in \eqref{eqn:probability-of-meeting-1}. To this end, for $\theta\in\mathbb{T}^d:=[-\pi,\pi]^d$, define
	\begin{equation}
		{F(\theta)}:=\sum_{j\in\Z^d}\e^{\ii(\theta,j)}\,b(0,j),\quad{\hat{F}(\theta)}
		:=\text{Re}({F(\theta)}),\quad {\tilde{F}(\theta)}:=\text{Im}({F(\theta)}).
	\end{equation}
	Then, for $j\in\Z^d$ and $t>0$,
	\begin{equation}
		\label{eqn:fourier-repr}
		\begin{aligned}
			\hat{a}_t(0,j) &= \frac{1}{(2\pi)^d}\int_{\mathbb{T}^d}\e^{-\ii (\theta,j)}\,\e^{-ct[1-{\hat{F}(\theta)}]}\,\d \theta,\\
			a_t(0,j) &= \frac{1}{(2\pi)^d}\int_{\mathbb{T}^d}\e^{-\ii (\theta,j)}\,\e^{-ct[1-{\hat{F}(\theta)}-\ii {\tilde{F}(\theta)}]}\,\d \theta.
		\end{aligned}
	\end{equation}
	Using that $a(i,0) = a(0,-i)$, $i\in\Z^d$, and inserting the above into \eqref{eqn:probability-of-meeting-1}, we obtain
	\begin{equation}
		\begin{aligned}
			\sum_{i\in\Z^d}a_{S(t)}(0,i)a_{S^\prime(t)}(0,i) 
			&= \frac{1}{(2\pi)^d}\int_{\mathbb{T}^d}\e^{-c[2M(t)+L(t)][1-{\hat{F}(\theta)}]}\cos(L(t){\tilde{F}(\theta)})\,\d\theta\\
			&=\frac{1}{(2\pi)^d}\int_{\mathbb{T}^d}\e^{-c[S(t)+S^\prime(t)][1-{\hat{F}(\theta)}]}\cos(L(t){\tilde{F}(\theta)})\,\d\theta\\
			&\leq \frac{1}{(2\pi)^d}\int_{\mathbb{T}^d}\e^{-c[S(t)+S^\prime(t)][1-{\hat{F}(\theta)}]}\,\d\theta\\
			&=\hat{a}_{S(t)+S^\prime(t)}(0,0),
		\end{aligned}
	\end{equation}
	where we use that $\tfrac{1}{(2\pi)^d}\sum_{j\in\Z^d}\e^{\ii(\theta-\theta^\prime,\,j)} = \delta(\theta-\theta^\prime)$, with $\delta(\cdot)$ the Dirac distribution (see e.g. \cite[Chapter 7]{Fol92}). Finally, combining the above with \eqref{eqn:average-intersection-time}--\eqref{eqn:probability-of-meeting}, we see that
	\begin{equation}
		I\leq \int_{0}^{\infty}\E^*_\eta\Big[\hat{a}_{S(t)+S^\prime(t)}(0,0)\Big]\,\d t
	\end{equation}
	and therefore it is necessary that
	\begin{equation}
		\label{eqn:nece-infin-inte}
		\int_{0}^{\infty}\E^*_\eta\Big[\hat{a}_{S(t)+S^\prime(t)}(0,0)\Big]\,\d t={\E^*_\eta\Big[\int_{0}^{\infty}\hat{a}_{S(t)+S^\prime(t)}(0,0)\,\d t\Big]}=\infty,
	\end{equation}
	{which by Remark~\ref{remark:recurrence-of-sym-kernel} is equivalent to
		\begin{equation}
			\int_{0}^{\infty}\hat{a}_t(0,0)\,\d t=\infty.
		\end{equation}
	}
	This proves the forward direction.
	
	To prove the converse, we first note that, because all the rates of absorption given by $(\tfrac{1}{N_i})_{i\in\Z^d}$ are such that \eqref{eqn:non-clumping criterion} holds and $\sup_{i\in\Z^d} K_i^{-1}<\infty$, whenever the two particles are on the same location, there is a positive probability of absorption that is uniformly bounded away from zero. Indeed, if $\nu^*(\eta)$ denote the absorption probability of $\xi^*$ when started from state $\eta$, by Theorem~\ref{thm:lower-bound-on-probability} we have that
	\begin{equation}
		\begin{aligned}
			&\inf_{i\in\Z^d}\nu^*([(i,1),(i,1)]) > 0,\\
			&\inf_{i\in\Z^d}\nu^*([(i,0),(i,1)]) = \inf_{i\in\Z^d}\nu^*([(i,0),(i,0)]) = \inf_{i\in\Z^d}\nu^*([(i,1),(i,0)]) >0,
		\end{aligned}
	\end{equation}
	where the last two equalities follow from a first-jump analysis of the process $\xi^*$ when started at the state $[(i,0),(i,0)],\ i\in\Z^d$. As a consequence, $\xi^*$ is absorbed with probability 1 if and only if, in the corresponding process where coalescence is switched off, the two particles infinitely often meet each other with probability 1. In other words, $\nu^*\equiv 1$ if and only if $I=\infty$, where $I$ is as in \eqref{eqn:average-intersection-time}, the average accumulated time spent by the two particles at the same location. However, by the symmetry of the kernel $a(\cdot\,,\cdot)$ and using Fubini's theorem, we have
	\begin{equation}
		I = \int_{0}^{\infty}\E^*_\eta\Big[a_{S(t)+S^\prime(t)}(0,0)\Big]\,\d t 
		= \int_{0}^{\infty}\E_\eta^*\Big[\hat{a}_{S(t)+S^\prime(t)}(0,0)\Big]\,\d t = {\E_\eta^*\Big[\int_{0}^{\infty}\hat{a}_{S(t)+S^\prime(t)}(0,0)\,\d t\Big]}
	\end{equation}
	and {thus, by Remark~\ref{remark:recurrence-of-sym-kernel}, if $\int_{0}^{\infty}\hat{a}_t(0,0)\,\d t=\infty$, then $I=\infty$}. This proves the backward direction.
\end{proof}

Now we are ready to prove Theorem~\ref{T.Coal_recc} with the help of Theorem~\ref{theorem: crit-for-abs-of-indep} and the results in Section~\ref{ss.concl}.

\begin{proof}[Proof of Theorem~\ref{T.Coal_recc}]
	Let $\nu(\eta)$ denote the absorption probability of the process $\xi$ (see Definition~\ref{definition: interacting RW1}) started at state $\eta\in G\times G$. Recall from Theorem~\ref{T.Dichotomy} and {Remark~\ref{remark:equivalence-of-absorption}} that the system clusters if and only if $\nu\equiv 1$. By the irreducibility of the process $\xi$, we have $\nu\equiv 1$ if and only if $\nu([(0,0),(0,0)])=1$. Now, since $\sup_{i\in\Z^d}K_i^{-1}<\infty$ and \eqref{eqn:non-clumping criterion} holds, we see that all the conditions of Theorem~\ref{theorem: necessity of coalescence of independent particles} are satisfied by virtue of Theorem~\ref{thm:lower-bound-on-probability}, and hence $\nu([(0,0),(0,0)])=1$ if and only if $\nu^*([(0,0),(0,0)])=1,$ where $\nu^*(\eta)$ denotes the absorption probability of the non-interacting two-particle process $\xi^*$ (see Definition~\ref{definition: independent RW}) started at state $\eta\in G\times G$. However, by the forward direction of Theorem~\ref{theorem: crit-for-abs-of-indep}, if $\nu^*([(0,0),(0,0)])=1$, then it is necessary that the symmetrised kernel $\hat{a}(\cdot\,,\cdot)$ is recurrent, and hence the forward direction is proved. Similarly, under the assumption of symmetry of the migration kernel, we can apply the converse direction of Theorem~\ref{theorem: necessity of coalescence of independent particles}, to conclude that if the transition kernel $a(\cdot\,,\cdot)$ (which is the same as the symmetrised transition kernel) is recurrent, then $\nu^*([(0,0),(0,0)])=1$, and so the backward direction follows as well. 
\end{proof}

{
	\begin{proof}[Proof of Corollary~\ref{C.dichotomy}]
		Recall from Remark~\ref{remark:moment-regularity} that the migration kernel $a(\cdot\,,\cdot)$ admits at least a $d$-th moment and is translation-invariant by assumption. Thus if $d>2$, then the kernel $\hat{a}(\cdot\,,\cdot)$, being symmetric by definition, is transient (by Polya's theorem), and hence clustering cannot take place by virtue of the forward direction of Theorem~\ref{T.Coal_recc}. Similarly, if $d \leq 2$ and $a(\cdot\,,\cdot)$ is symmetric, then $a(\cdot\,,\cdot)$ is recurrent, and so the claim follows from the backward direction of Theorem~\ref{T.Coal_recc}.
	\end{proof}
}


{
	\section{Discussion}
}

{Stochastic models describing genetic evolution of \emph{finite} populations under various evolutionary forces remain a challenge in population genetics. The presence of a seed-bank can complicate the analysis even further. In recent years, stochastic duality has proven to be a very useful mathematical tool, particularly in the field of interacting particle system, for tackling technical complications and doing explicit computations. On the one hand, we aim to create a bridge between interacting particle system and mathematical population genetics by including dormancy into existing well-known particle systems. On the other hand, we hope to combine this approach with the recently developed theory of duality to reveal delicate structures and related interesting properties of the interacting particle system that lie hidden and are often lost in the process of taking the large-colony-size limit.}

{In \cite{HN01}, we heavily rely on duality to prove our results for to the process $Z$. In a subdivided population, the ancestral dual process in the presence of resampling and migration is generally described by the \emph{structured coalescent process}. This process, which is by now well-understood, was originally derived as the genealogical process in the context of geographically structured \emph{large} populations under Wright-type reproduction and migration (see e.g., \cite{H94,W2001} and \cite{S15}). Even though lineages move independently in the structured coalescent, the genealogies of a sample taken from subdivided and finite populations with constant size are correlated \cite{N90,H94}. These correlations arise due to the imposition of \emph{finite} and \emph{constant} (in time) population sizes, and vanish when the large-population-size limit is taken.}

{As can be seen in \cite[Definition 3.7]{HN01}, the ancestral dual process $Z^*$ is no exception, and lineages in the dual indeed show a \emph{repulsive} interaction. Due to the incorporation of dormancy, lineages can also adopt one of two states: active and dormant. The presence of these correlations and of dormant periods in the lineages make the dual process $Z^*$ interesting but tricky to analyse. Consequently, in the present paper we take a different route to address the dichotomy of coexistence versus clustering. More precisely, instead of directly exploiting the clustering criterion given in terms of the original two-particle dual process (equivalently, the process $\xi$ in Definition~\ref{definition: interacting RW1}), we find an alternative clustering criterion that is relatively easy to deal with. We achieve this by comparing the original two-particle dual $\xi$ with two \emph{auxiliary} two-particle duals processes $\hat{\xi}$ and $\xi^*$ (see Definition~\ref{definition: interacting RW2} and Definition~\ref{definition: independent RW}), which are simplified versions of $\xi$. In particular, we obtain $\hat{\xi}$ from $\xi$ by switching off the repulsive interaction present in the migration mechanism of an active particle and removing the coalescence of active particles from different locations, while $\hat{\xi}$ is further simplified to $\xi^*$, the \emph{independent RW} process, by turning off the only interaction that takes place between an active and a dormant particle located at the same position. The comparison technique employed in Section~\ref{s.duals} to estimate the absorption probabilities for $\xi,\hat{\xi},\xi^*$ is similar to that in \cite{G14}, where a connection is made between infinitesimal generators of the Wright-Fisher diffusion and the $\Lambda$-Fleming-Viot process, based on methods involving Lyapunov functions to characterise fixation probabilities. Similar techniques are used in the literature of interacting particle systems to derive correlation inequalities and related properties (see e.g., \cite{GRV10}). It is worth emphasising that our results are valid for any choice of the sizes $(N_i)_{i\in\Z^d}$ and $(M_i)_{i\in\Z^d}$ of active and dormant populations, subject to the mild criteria we imposed. Such generalities are rare and suggest that other problems can perhaps be approached in a similar way.
}
\bigskip
\appendix
\gdef\thesection{\Alph{section}} 
\makeatletter
\renewcommand\@seccntformat[1]{\csname the#1\endcsname.\hspace{0.5em}}
\makeatother
\section{Two-particle dual and alternative representation}
\label{aps-basic-dual}
\noindent
{In this appendix, we give a short description of the original dual process $\tilde{Z}$ started with two particles, which was introduced in full generality as a configuration process $Z^*$ in \cite[Section 3.2]{HN01}. Further, we briefly outline the derivation of the \emph{interacting RW1} process $\xi$ defined in Definition~\ref{definition: interacting RW1} from the configuration process $\tilde{Z}$, and show that the absorption of $\xi$ and coalescence of the two particles in $\tilde{Z}$ are basically equivalent.} 
\begin{definition}{\bf [Two-particle dual]}
	{\rm \label{definition: two-particles dual}
		The two-particle dual process 
		\begin{equation}
			\tilde{Z} := (\tilde{Z}(t))_{t \geq 0}, \qquad \tilde{Z}(t) := (\tilde{n_i}(t),\tilde{m_i}(t))_{i\in\Z^d},
		\end{equation} 
		is the continuous-time Markov chain with state space 
		\begin{equation}
			\tilde{\X} := \Big\{(\tilde{n_i},\tilde{m_i})_{i\in\Z^d} \in 
			\prod_{i\in\Z^d}[N_i]\times[M_i]\colon\,\sum_{i\in\Z^d}(\tilde{n_i}+\tilde{m_i})\leq 2\Big\}
		\end{equation}
		and with transition rates
		\begin{equation}
			\label{rates: dual}
			\begin{aligned}
				&(n_k,m_k)_{k\in\Z^d} \to\\
				&\begin{cases}
					\displaystyle
					(n_k,m_k)_{k\in\Z^d} - \vec{\delta}_{i,A} ,\ \ \,\,\quad\quad\quad\text{ at rate } 
					\tfrac{2a(i,i)}{N_i} \binom{n_i}{2}\mathbf{1}_{\{n_i\geq 2\}}
					+\sum_{j\in\Z^d\setminus\{i\}} \tfrac{n_ia(i,j)n_j}{N_j} \text{ for } i\in\Z^d,\\
					(n_k,m_k)_{k\in\Z^d} - \vec{\delta}_{i,A}+\vec{\delta}_{i,D}, \quad \text{ at rate } 
					\tfrac{\lambda n_i(M_i-m_i)}{M_i} \ \qquad \text{ for } i\in\Z^d,\\
					(n_k,m_k)_{k\in\Z^d} + \vec{\delta}_{i,A}-\vec{\delta}_{i,D}, \quad\text{ at rate } 
					\tfrac{\lambda (N_i-n_i)m_i}{M_i}\ \qquad \,\text{ for } i\in\Z^d,\\
					(n_k,m_k)_{k\in\Z^d} - \vec{\delta}_{i,A} +\vec{\delta}_{j,A}, \quad\text{ at rate } 
					\tfrac{n_ia(i,j)(N_j-n_j)}{N_j} \quad \,\text{ for } i\neq j\in\Z^d,
				\end{cases}
			\end{aligned}
		\end{equation}
		where for $i\in\Z^d$ {the configurations $\vec{\delta}_{i,A}$, $\vec{\delta}_{i,D}$ are defined as}
		\begin{equation}
			\label{kron_del}
			{\vec{\delta}_{i,A}:= (\mathbf{1}_{\{n=i\}},0)_{n\in\Z^d},\quad
				\vec{\delta}_{i,D}:=(0,\mathbf{1}_{\{n=i\}})_{n\in\Z^d},}
		\end{equation}
		and for two configurations $\eta_1 = (\bar{X}_i,\bar{Y}_i)_{i\in\Z^d}$ and $\eta_2 = (\hat{X}_i,\hat{Y}_i)_{i\in\Z^d}$, $\eta_1\pm\eta_2:=(X_i,Y_i)_{i\in\Z^d}$ is defined componentwise by
		\begin{equation}
			\label{eq:addition-subtraction}
			\begin{aligned}
				X_i &= (\bar{X_i}\,\pm\,\hat{X_i})\mathbf{1}_{\{0\leq \bar{X}_i\pm\hat{X}_i\leq N_i\}}
				+N_i\,\mathbf{1}_{\{ \bar{X}_i\pm\hat{X}_i> N_i\}},\\
				Y_i &= (\bar{Y_i}\,\pm\,\hat{Y}_i)\mathbf{1}_{\{0\leq \bar{Y}_i\pm\hat{Y}_i\leq M_i\}}
				+M_i\,\mathbf{1}_{\{ \bar{Y}_i\pm\hat{Y}_i> M_i\}}.
			\end{aligned}
		\end{equation}
		The support of the distribution of $\tilde{Z}(0)$ is contained in 
		\begin{equation}
			\tilde{\X}_0 := \Big\{(\tilde{n_i},\tilde{m_i})_{i\in\Z^d} \in 
			\prod_{i\in\Z^d}[N_i]\times[M_i]\colon\,\sum_{i\in\Z^d}(\tilde{n_i}+\tilde{m_i}) =  2\Big\}.
		\end{equation}
	} \hfill $\Box$
\end{definition}

\noindent
Here, $\tilde{n_i}(t)$ and $\tilde{m_i}(t)$ are the number of active and dormant particles at site $i\in\Z^d$ at time $t$. The first transition describes the coalescence of an active particle at site $i$ with active particles at other sites. The second and third transition describe the switching between the active and the dormant state of the particles at site $i$. The fourth transition describes the migration of an active particle from site $i$ to site $j$.

Let $\tilde{\X}_1$ be the set of configurations containing a single particle, i.e.,
\begin{equation}
	\tilde{\X}_1:=\Big\{(\tilde{n_i},\tilde{m_i})_{i\in\Z^d} \in 
	\tilde{\mathcal{X}}\colon\,\sum_{i\in\Z^d}(\tilde{n_i}+\tilde{m_i}) = 1\,\Big\},
\end{equation}
and let $\tilde{\tau}$ be the first time at which coalescence has occurred, i.e.,
\begin{equation}
	\label{definition: coalescence}
	\tilde{\tau}=\inf\{t\geq0\colon\,(\tilde{n_i}(t),\tilde{m_i}(t))_{i\in\Z^d} \in \tilde{\X}_1\}.
\end{equation}
As indicated earlier in Section~\ref{ss.basicduals}, we are only required to analyse the coalescence probability of two dual particles and thus, it boils down to lumping all the configurations in $\tilde{\X}_1$ into a single state $\circledast$ and consider the resulting lumped process. Note that, on the event $\{\tilde{\tau}< s\}$, the process $(\tilde{Z}(t))_{t\geq s}$ a.s.\ stays in $\tilde{\X}_1$. Therefore the lumped process is a well-defined continuous-time Markov chain with state space $\tilde{\mathcal{X}}_0\cup\{\circledast\}$, where $\circledast$ is an absorbing state.

With a little abuse of notation, from here onwards we denote the lumped process by $(\tilde{Z}(t))_{t\geq 0}$. We give the formal description of this process in a definition.

\begin{definition}{\bf [Lumped two-particle dual]}
	\label{definition: lumped dual}
	{\rm The lumped two-particle dual process 
		\begin{equation}
			\tilde{Z} := (\tilde{Z}(t))_{t \geq 0}
		\end{equation} 
		is the continuous-time Markov chain with state space 
		\begin{equation}
			\tilde{\X} := \Big\{(n_i,m_i)_{i\in\Z^d} \in 
			\prod_{i\in\Z^d}[N_i]\times[M_i]\colon\,\sum_{i\in\Z^d}(n_i+m_i)= 2\Big\}\bigcup\Big\{\circledast\Big\}
		\end{equation}
		and with transition rates
		\begin{equation}
			\label{rates: lumped-dual}
			\begin{aligned}
				&(n_k,m_k)_{k\in\Z^d} \to\\
				&\begin{cases}
					\displaystyle
					\circledast, &\text{ at rate } \displaystyle
					\sum_{i\in\Z^d}\Big[\tfrac{2a(0,0)}{N_i}\binom{n_i}{2}\mathbf{1}_{\{n_i\geq 2\}}
					+\sum_{j\in\Z^d\setminus\{i\}}\tfrac{n_ia(i,j)n_j}{N_j}\Big] \text{ for } i\in\Z^d,\\
					(n_k,m_k)_{k\in\Z^d} - \vec{\delta}_{i,A}+\vec{\delta}_{i,D}, &\text{ at rate } 
					\tfrac{\lambda n_i(M_i-m_i)}{M_i} \ \qquad \text{ for } i\in\Z^d,\\
					(n_k,m_k)_{k\in\Z^d} + \vec{\delta}_{i,A}-\vec{\delta}_{i,D}, &\text{ at rate } 
					\tfrac{\lambda (N_i-n_i)m_i}{M_i}\ \qquad \,\text{ for } i\in\Z^d,\\
					(n_k,m_k)_{k\in\Z^d} - \vec{\delta}_{i,A} +\vec{\delta}_{j,A}, &\text{ at rate } 
					\tfrac{n_ia(i,j)(N_j-n_j)}{N_j} \quad \,\text{ for } i\neq j\in\Z^d,
				\end{cases}
			\end{aligned}
		\end{equation}
		where, for $i\in\Z^d,$ $\vec{\delta}_{i,A}$ and $\vec{\delta}_{i,D}$ are as in \eqref{kron_del}.
	} \hfill $\Box$
\end{definition}

\noindent
We write $\tilde{\P}^\eta$ to denote the law of the process $\tilde{Z}$ started from $\eta\in\mathcal{X}$. Note that, by construction, the coalescence time $\tilde{\tau}$ is now same as the absorption time of the process $\tilde{Z}$. In the following proposition, we show that the configuration process $\tilde{Z}$ is an alternative representation of the coordinate process $\xi$ defined in Definition~\ref{definition: interacting RW1}.


\begin{proposition}{\bf [Equivalence between $\tilde{Z}$ and $\xi$]}
	\label{Prop: equivalence between lumped dual and interacting RW1}
	Let $\xi = (\xi(t))_{t\geq 0}$ be the process defined in Definition \ref{definition: interacting RW1} with initial distribution $\mu$. Let $\phi\colon\,\mathcal{S}\to\tilde{\mathcal{X}}$ be the map defined by
	\begin{equation}
		\phi(\eta) := 
		\begin{cases}
			(\alpha\delta_{k,i}+\beta\delta_{k,j},(1-\alpha)\delta_{k,i}+(1-\beta)\delta_{k,j})_{k\in\Z^d}, 
			&\text{ if }\eta=[(i,\alpha),(j,\beta)]\neq\circledast,\\
			\circledast, 
			& \text{ otherwise}.
		\end{cases}
	\end{equation}
	For $t\geq 0$, let $\tilde{Z}(t):=\phi(\xi(t))$. Then the process $(\tilde{Z}(t))_{t\geq 0}$ is the lumped dual process defined in Definition \ref{definition: lumped dual}, and its initial distribution is the push-forward of $\mu$ under the map $\phi$. {Furthermore, $\tilde{Z}$ is absorbed to $\circledast$ if and only if $\xi$ is}.
\end{proposition}

\begin{proof}
	Due to the Assumption~\ref{assump:non-trivial colony-sizes}, we see that $\phi(\eta)\in\tilde{\mathcal{X}}$, and so $\tilde{Z}(t)\in\X$ for all $t\geq 0$, and $\phi$ is onto. For $\eta\in\S$, define 
	\begin{equation}
		\bar{\eta} := 
		\begin{cases}
			[(j,\beta),(i,\alpha)], & \text{ if }\eta=[(i,\alpha),(j,\beta)]\neq\circledast,\\
			\circledast, & \text{ otherwise}.
		\end{cases}
	\end{equation}
	Note that $\phi^{-1}(\phi(\eta)) = \{\eta,\bar{\eta}\}$. Let $Q(\eta_1,\eta_2)$ denote the transition rate from $\eta_1$ to $\eta_2$ for the process $\xi$, where $\eta_1\neq \eta_2\in\S$. Furthermore, let $z_1\neq z_2\in\tilde{\X}$ be fixed and $\eta_1\in\S$ be such that $\phi(\eta_1)=z_1$. Since $Q(\eta_1,\eta_2)=Q(\bar{\eta}_1,\bar{\eta}_2)$ for any $\eta_1\neq\eta_2\in\S$, we have
	\begin{equation}
		\label{eqn: dynkin-criterion}
		\begin{aligned}
			\sum_{\eta\in\phi^{-1}(z_2)} Q(\eta_1,\eta) = \sum_{\eta\in\phi^{-1}(z_2)} Q(\bar{\eta}_1,\eta).
		\end{aligned}
	\end{equation}
	Hence the Dynkin criterion for lumpablity is satisfied, and $\phi$ preserves the Markov property. So $\tilde{Z}$ is a Markov process on $\tilde{X}$. We can easily verify that the sum in \eqref{eqn: dynkin-criterion} is indeed the transition rate from $z_1$ to $z_2$ defined in \eqref{rates: lumped-dual}. Thus, $\tilde{Z}$ is the lumped dual process defined in Definition~\ref{definition: lumped dual}. Clearly, the distribution of $\tilde{Z}(0)$ is $\mu\circ \phi^{-1}$. {The second claim trivially follows, since $\phi(\eta)=\circledast$ if and only if $\eta=\circledast$}.
\end{proof}

\bibliographystyle{elsarticle-num}

 
\end{document}